\documentclass[journal,twoside,web]{ieeecolor}
\usepackage{generic}
\usepackage{cite}
\usepackage{amsmath,amssymb,amsfonts}
\usepackage{algorithmic}
\usepackage{graphicx}
\usepackage{textcomp}
\usepackage{xcolor}
\usepackage{booktabs}
\usepackage{mathtools}
\usepackage{tabu}
\usepackage{tabularx}
\usepackage{array}
\usepackage{float}
\usepackage{subcaption}
\usepackage{upgreek}
\usepackage{tikz}
\usetikzlibrary{positioning}
\usepackage{multirow}
\usepackage{booktabs}
\usepackage{lipsum}
\usepackage{cite}

\newcolumntype{C}{>{\centering\let\newline\\\arraybackslash\hspace{0pt}}m{2.3cm}}

\newtheorem{theorem}{Theorem}
\newtheorem{lemma}{Lemma}
\newtheorem{proposition}{Proposition}
\newtheorem{problem}{Problem}
\newtheorem{remark}{Remark}

\newtheorem{prooftheorem}{Proof of Theorem}
\stepcounter{prooftheorem}
\newtheorem{prooflemma}{Proof of Lemma}
\newtheorem{proofproposition}{Proof of Proposition}

\def\BibTeX{{\rm B\kern-.05em{\sc i\kern-.025em b}\kern-.08em
    T\kern-.1667em\lower.7ex\hbox{E}\kern-.125emX}}
\makeatletter
\let\NAT@parse\undefined
\makeatother
\usepackage{hyperref}[backref]
\hypersetup{
    colorlinks,
    linkcolor={blue!70!black},
    citecolor={red!70!black},
    urlcolor={blue!80!black}
}

\begin{document}

\title{New Dualities in Linear Systems and Optimal Output Control under Bounded Disturbances}

\author{
\thanks{{\color{lightgray} Manuscript info will be here. Manuscript info will be here. Manuscript info will be here. Manuscript info will be here. Manuscript info will be here. Manuscript info will be here. Manuscript info will be here. Manuscript info will be here. Manuscript info will be here. Manuscript info will be here. Manuscript info will be here. Manuscript info will be here. Manuscript info will be here. Manuscript info will be here. Manuscript info will be here.
Manuscript info will be here. Manuscript info will be here.
Manuscript info will be here. Manuscript info will be here.
Manuscript info will be here. Manuscript info will be here.
Manuscript info will be here. Manuscript info will be here.
Manuscript info will be here. }}
Alexey Peregudin,
\thanks{A. Peregudin is with the ITMO University, 49 Kronverkskiy ave, Saint Petersburg, 197101, and the Institute for Problems of Mechanical Engineering Russian Academy of Sciences (IPME RAS), 61 Bolshoy Prospekt V.O., St.-Petersburg, 199178, Russia (e-mail: peregudin@itmo.ru).}
Igor Furtat
\thanks{I. Furtat is with the Institute for Problems of Mechanical Engineering Russian Academy of Sciences (IPME RAS), 61 Bolshoy Prospekt V.O., St.-Petersburg, 199178, Russia (e-mail: cainenash@mail.ru).}
\thanks{{\color{lightgray}  }}
}

\maketitle

\begin{abstract}
In this paper, we introduce novel equations that are dual to the ones of the well-known invariant ellipsoids method. These equations yield ellipsoids with newly established geometrical interpretations and connections to linear system norms. The established duality leads to the optimal synthesis results for state-feedback control, filtering, and output-feedback control problems in the presence of bounded disturbances. The proposed output-feedback control solution is demonstrated to be optimal and surpass prior sub-optimal results.

\end{abstract}

\begin{IEEEkeywords}
 Attractive ellipsoids, invariant ellipsoids, linear systems, optimal control, output feedback, bounded disturbances, system norms.
\end{IEEEkeywords}

\section{Introduction}
\label{sec:introduction}

\IEEEPARstart{A}TTENUATION of bounded disturbances is a major challenge in control systems. Such disturbances are 
only assumed to be non-stochastic, but not necessarily decaying, generated by a linear system, or of finite $L_2$-norm.  Optimal control strategies that are widely used for other control problems (such as LQR/$\mathcal{H}_2$ and $\mathcal{H}_\infty$-control) are not suitable for bounded deterministic disturbances as they minimize different cost functions not relevant to the task.

The invariant (attractive) ellipsoids method is a way to analyze and design control laws for systems subject to bounded disturbances of $L_\infty$-class. The method aims to minimize the set of states reachable under these disturbances by minimizing the size of the ellipsoid that approximates it. The fundamentals of the invariant ellipsoids method are first outlined in \cite{b12-boyd}, \cite{b5-abedor}, and later developed in \cite{b1-nazin, b13-khlebnikov_filter, b2-topunov}. The method is used for state-feedback control \cite{b1-nazin}, observer design \cite{b13-khlebnikov_filter}, and output-feedback control \cite{b2-topunov}. The further development of the method was associated with its adaptation to more complex tasks such as control of some classes of nonlinear systems \cite{e1}, network systems control \cite{x2}, adaptive control \cite{y1}, robust control \cite{x6}, \cite{x7}, and sliding mode control \cite{e3}, \cite{y3}. The classical monograph \cite{b3-poznyak} on attractive ellipsoids summarize some of these results as well as their extensions. 

Among recent works, it is worth mentioning \cite{y6, e5, e6, y4, y5, y7}. In \cite{y6} attractive ellipsoids are applied to the simultaneous localization and mapping problem, in \cite{y4} and \cite{y5} the method is used for quadrotor control, and in \cite{y7} it is applied for controlling time-varying polytopic systems.

Despite its popularity, there are some issues with the invariant/attractive ellipsoids method. {{}Firstly, while the geometric interpretation of the invariant ellipsoids equation is well-known, the same cannot be said for its dual counterpart, which, to the best of the authors' knowledge, has never been established or used before.} Secondly, the method relies heavily on optimization, resulting in controller synthesis procedures that are never exact and always given as optimization problems with LMI constraints. Finally, the output-feedback controller in \cite{b2-topunov}, \cite{b3-poznyak} derived from this method is known to be sub-optimal, while the optimal solution to the output-feedback control problem was never proposed.

The main contributions of this paper are as follows:

\begin{itemize}
    \item the interpretation of {{}equation \eqref{Q1} in Section \ref{S2B}}, which is dual to the invariant ellipsoids equation \eqref{Pinf}, is provided and the relationship between the solutions of these equations and system norms is established;
    \item based on the established {{}duality relations}, exact equations for the optimal state-feedback controller and filter are presented, {{}which have fewer variables compared to previously known LMI-based methods, therefore enabling faster and more precise computations;}
    \item for the first time, the optimal solution to the problem of output-feedback control under bounded disturbances is proposed.
\end{itemize}

This paper is organized as follows. Section \ref{S2} focuses on analyzing system performance in terms of reachability and observability and establishes results that are dual and symmetric to the known ones. In Section \ref{S3}, the state-feedback and filtering problems are addressed using the perspective described in Section \ref{S2}. Section \ref{S4} introduces a novel solution to the output-feedback control problem. Section \ref{S5} compares new results with the previous ones both theoretically and numerically. Section \ref{S6} discusses possible improvements to the method, and Section \ref{S7} provides the conclusions. The proofs of all propositions, theorems, as well as technical lemmas can be found in the \hyperlink{appendix}{Appendix}.

\subsection{Notation}
$\mathbb{R}$ is the set of all real numbers, $\mathbb{R}^{m \times n}$ is the set of all $m\times n$ matrices with real entries, $\mathbb{R}^n \vcentcolon = \mathbb{R}^{n \times 1}$. The transpose of $A \in \mathbb{R}^{m \times n}$ is denoted by $A^\top$. We use the Euclidean norm $|v| \vcentcolon = \sqrt{v^\top v}$ for $v \in \mathbb{R}^n$. Matrix $A \in \mathbb{R}^{n \times n}$ is said to be stable iff all its eigenvalues have strictly negative real parts. $\lambda_{\text{max}} (A)$ stands for the maximum eigenvalue of a symmetric matrix $A$, $\sigma_{\text{max}} (A)$  stands for the maximum singular value of a general matrix $A$. The ordering symbols $\succ$, $\prec$, $\succeq$, $\preceq$ are used in the sense of matrix definiteness, e.g. $A \succ B$ means that both $A$ and $B$ are symmetric and $A-B$ is positive definite.

We define the set of {{}$n$-dimensional} signals as
\begin{equation*}
    \mathcal{F}^n(T) \vcentcolon = \left\{ f : [0, T] \to \mathbb{R}^n, \;\: \,  f \text{ is measurable}\right\}, \\ 
\end{equation*}
and for $f \in \mathcal{F}^n(T)$ and $r\ge 1$ we use the norms
\begin{equation*}
\|f\|_r \vcentcolon = \left( \int_{0}^{T} |f(t)|^r dt \right)^{1/r}\!\!\!, \quad \|f\|_\infty \vcentcolon = \underset{t \in [0, T]} {\operatorname{ess \, sup}} \, |f(t)|, 
\end{equation*}
when the corresponding values are well-defined. Note that these norms are only applied to finite-time signals, but the systems will be studied in infinite-time. The transition is achieved through a set-theoretical limiting process (see the definitions of $\mathcal{R}_p$ and $\mathcal{O}_q$ in Section \ref{S2}).

\section{Analysis: Easily Reachable and Hardly Observable Sets, System Norms}
\label{S2}

\subsection{Preliminaries on Reachability and Observability}

Consider a linear time-invariant strictly proper system
\begin{equation*}
    \mathcal{S} : \; \left\{ 
    \begin{aligned}
        &\dot x = Ax + Bu, \\
        &y = Cx,
    \end{aligned}
    \right.
\end{equation*}
where {{} $x(t) \in \mathbb{R}^n$, $u(t) \in \mathbb{R}^m$, $y(t) \in \mathbb{R}^k$, and $A$, $B$, $C$ are real matrices of corresponding sizes. Also consider} its Input-to-State and State-to-Output components
\begin{equation*}
    \mathcal{S}_u : \; \dot x = Ax + Bu, \quad \quad \mathcal{S}_y : \; 
        \dot x = Ax, \;
        y = Cx.
\end{equation*}
Assume that $B$ has full column rank and $C$ has full row rank. This can always be achieved by removing linearly dependent inputs and outputs from the model.

Define $\mathcal{R}(T)$ as the set of all states that are reachable at a given time $T$ with unconstrained input, i.e.
\begin{equation*}
    \mathcal{R}(T) \vcentcolon = \left\{ x(T) \in \mathbb{R}^n  \mid \mathcal{S}_u, \; x(0) = 0, \; u \in \mathcal{F}^m(T)\right\}. 
\end{equation*}
Note that for $T>0$ all such sets coincide, so we will denote them by $\mathcal{R}$. In fact, $\mathcal{R}$ is the \textit{reachable subspace} for $\mathcal{S}$. 

Define $\mathcal{O}(T)$ as the set of all initial states that result in an identically zero output up to a given time $T$, i.e.
\begin{equation*}
    \mathcal{O}(T) \vcentcolon = \left\{ x(0)\in \mathbb{R}^n  \mid \mathcal{S}_y, \; y(t) \equiv 0, \; y \in \mathcal{F}^k(T)\right\}. 
\end{equation*}
Again, for $T>0$ all such sets coincide, so we will denote them by $\mathcal{O}$. In fact, $\mathcal{O}$ is the \textit{unobservable subspace} for $\mathcal{S}$. 

For $p\ge 1$ define $\mathcal{R}_p(T)$ as the set of all states reachable at a time $T$ with an input with no more than a unit $p$-norm, i.e.
\begin{equation*}
    \mathcal{R}_p(T) \vcentcolon = \left\{ x(T) \mid\mathcal{S}_u, \; x(0) = 0, \; u \in \mathcal{F}^m(T), \; \|u\|_p \le 1 \right\}.
\end{equation*}
It is straightforward to show that $\mathcal{R}_p(T)$ is convex and is strictly convex if the system is completely controllable. The set is expanding, i.e. for $T_2 > T_1$ one has $\mathcal{R}_p(T_1) \subset \mathcal{R}_p(T_2)$. Define the \textit{{{}easily reachable} set} (with respect to $p$-norm) for a system $\mathcal{S}$ as
\begin{equation*}
    \mathcal{R}_p \vcentcolon = \overline{\bigcup_{T \ge 0}\mathcal{R}_p(T)},
\end{equation*}
where the overline represents the topological closure. We have the following inclusion
\begin{equation*}
    \mathcal{R}_p(T) \subset \mathcal{R}_p \subset \mathcal{R}.
\end{equation*}

Note that if $\mathcal{S}$ is completely controllable, then $\mathcal{R}_p$ is bounded iff the matrix $A$ is stable.

Analogously, for $q\ge 1$ define $\mathcal{O}_q(T)$ as the set of all initial states that result in an output with no more than a unit $q$-norm up to a time $T$, i.e.
\begin{equation*}
     \mathcal{O}_q(T) \vcentcolon = \left\{ x(0) \mid \mathcal{S}_y, \;  y \in \mathcal{F}^k(T), \; \|y\|_q \le 1\right\}. 
\end{equation*}
The set $\mathcal{O}_q(T)$ is convex and is strictly convex, if the system is completely observable. The set is contracting, i.e. for $T_2 > T_1$ one has $\mathcal{O}_q(T_2) \subset \mathcal{O}_q(T_1)$. Define the \textit{{{}hardly observable} set} (with respect to $q$-norm)  for a system $\mathcal{S}$ as
\begin{equation*}
    \mathcal{O}_q \vcentcolon = \bigcap_{T \ge 0}\mathcal{O}_q(T).
\end{equation*}
Then we get the following inclusion
\begin{equation*}
    \mathcal{O} \subset \mathcal{O}_q \subset \mathcal{O}_q(T).
\end{equation*}
Note that if $\mathcal{S}$ is completely observable, then $0$ is an interior point in $\mathcal{O}_q$ iff the matrix $A$ is stable.

\subsection{Ellipsoidal Approximations of Easily Reachable and Hardly Observable Sets} \label{S2B}

From now on we assume that $\mathcal{S}$ is completely controllable, completely observable and that the matrix $A$ is stable. Then it is well-known (see \cite{b9-optimalbook}, \cite{b16-usual_ellipsoids}) that 
\begin{equation*}
    \mathcal{R}_2 = \left\{ x \mid x^\top P^{-1} x \le 1 \right\}, \quad
    \mathcal{O}_2 = \left\{ x \mid x^\top Q x \le 1 \right\},
\end{equation*}
where {\color{black}$P$, $Q \succ 0$ are the controllability and observability Gramians, i.e. the unique solutions to
\begin{equation}
\label{P2Q2}
    AP+PA^\top+BB^\top = 0,  \quad 
    QA+A^\top Q+C^\top C = 0.
\end{equation}}
It means that for the case $p=q=2$ both {{}easily reachable} and {{}hardly observable} sets are exact ellipsoids. 

Ellipsoidal approximations for the set $\mathcal{R}_\infty$ are the subject of the invariant ellipsoids method \cite{b5-abedor}, \cite{b1-nazin}, \cite{b3-poznyak}. The following theorem is known.

\begin{theorem}[\hspace{1sp}\cite{b5-abedor}, \cite{b1-nazin}]
If $\alpha >0$ and $P_\alpha \succ 0$ are such that
\begin{equation} 
\label{Pinf}
    AP_\alpha+P_\alpha A^\top+\alpha P_\alpha + \frac{1}{\alpha}BB^\top = 0,
\end{equation}
then we have the inclusion
\begin{equation*}
    \mathcal{R}_\infty \subset  \left\{ x \mid x^\top P_\alpha^{-1} x \le 1 \right\}.
\end{equation*}
\label{thm:Pinf}
\vspace{-\baselineskip}
\end{theorem}

In the literature, the sets $\left\{ x \mid x^\top P_\alpha^{-1} x \le 1 \right\}$ are commonly referred to as invariant ellipsoids \cite{b1-nazin}, \cite{b2-topunov} or  attractive ellipsoids \cite{b3-poznyak}. However, the dual equation \eqref{Q1} has not been explored before and its geometric meaning was unknown.  It is for the first time that the duality between $\mathcal{R}_\infty$ and $\mathcal{O}_1$ is established through the following theorem. 

\begin{theorem} 
\label{thm:Q1}
If $\alpha >0$ and $Q_\alpha \succ 0$ are such that
\begin{equation} 
\label{Q1}
    Q_\alpha A+A^\top Q_\alpha +\alpha Q_\alpha + \frac{1}{\alpha}C^\top C = 0,
\end{equation}
then we have the inclusion
\begin{equation*}
     \left\{ x \mid x^\top Q_\alpha x \le 1 \right\} \subset \mathcal{O}_1.
\end{equation*}
\vspace{-\baselineskip}
\end{theorem} 

The proof of Theorem \ref{thm:Q1}, along with the proofs of other theorems and propositions, can be found in the \hyperlink{appendix}{Appendix}.

We have thus established the dual relationship between the ellipsoidal approximations of $\mathcal{R}_\infty$ and $\mathcal{O}_1$. To provide a complete picture (although it is not needed for the control design in later sections), we will demonstrate the dual relationship between the ellipsoidal approximations of $\mathcal{R}_1$ and $\mathcal{O}_\infty$.

In \cite{b12-boyd}, an ellipsoidal approximation for $\mathcal{O}_\infty$ is given without proof. We state it as a theorem and provide proof, as well as its counterpart for $\mathcal{R}_1$, which has not been previously presented.

\begin{theorem}[\hspace{1sp}\cite{b12-boyd}]
\label{thm:Qinf}
    If $\tilde Q \succ 0$ is such that 
    \begin{equation}
    \label{Qinf}
       \tilde QA+A^\top \tilde Q \prec 0, \quad  \tilde Q \succeq  C^\top C, 
    \end{equation}
    then we have the inclusion 
    \begin{equation*}
        \left\{ x \mid x^\top \tilde Q x \le 1 \right\} \subset \mathcal{O}_\infty.
    \end{equation*}
    \vspace{-\baselineskip}
\end{theorem}


\begin{theorem}
\label{thm:P1}
    If $\tilde P \succ 0$ is such that 
    \begin{equation}
    \label{P1}
        A\tilde P+\tilde PA^\top \prec 0, \quad  \tilde P \succeq  B B^\top, 
    \end{equation}
    then we have the inclusion 
    \begin{equation*}
        \mathcal{R}_1 \subset  \left\{ x \mid x^\top \tilde P^{-1} x \le 1 \right\} .
    \end{equation*}
    \vspace{-\baselineskip}
\end{theorem}

\begin{table*}[bt]
\caption{Comparison table. We use the symbol {\color{subsectioncolor} ${\color{subsectioncolor} \star}$}  to indicate new concepts and results that are the contribution of this paper.}
\centering
    \begin{tabular}{llll}
    \toprule
          Reachability equations &  
          $AP_\alpha+P_\alpha A^\top+\alpha P_\alpha + \frac{1}{\alpha}BB^\top = 0$ & 
          $AP+PA^\top+BB^\top = 0$ &  
          $A\tilde P+\tilde PA^\top \prec 0, \; \; \tilde P \succeq  B B^\top$ \\[0.15cm] 
          Observability equations & 
          $Q_\alpha A+A^\top Q_\alpha +\alpha Q_\alpha + \frac{1}{\alpha}C^\top C = 0$ & 
          $QA+A^\top Q+C^\top C= 0$ &  
          $\tilde Q A+A^\top \tilde Q \prec 0, \; \;  \tilde Q \succeq C^\top C$ \\  
          \midrule
        {{}Easily reachable} set approx. & 
         $\mathcal{R}_\infty \subset  \left\{ x \mid x^\top P_\alpha^{-1} x \le 1 \right\}$ & 
         $\mathcal{R}_2 = \left\{ x \mid x^\top P^{-1} x \le 1 \right\}$ & 
         $\mathcal{R}_1 \subset  \left\{ x \mid x^\top \tilde P^{-1} x \le 1 \right\} \; {\color{subsectioncolor} {\color{subsectioncolor} \star}}$ \\[0.15cm] 
          {{}Hardly observable} set approx. & 
         $ \quad \quad \: \, {\color{subsectioncolor} \star} \; \left\{ \, x \, \mid 
 \, x^\top Q_\alpha x \le 1  \right\} \subset \mathcal{O}_1   $ & 
         $\quad \quad \; \; \:  \left\{\;  x \; \mid \; x^\top Q x \le 1 \ \right\} = \mathcal{O}_2 $ & 
         $  \quad \quad \; \; \: \left\{ \; x \; \mid \; x^\top \tilde Q x \le 1 \ \right\} \subset \mathcal{O}_\infty $ \\
         \midrule 
        Largest eigenvalue norms &
          $\displaystyle \| \mathcal{S}\|_\ast^2  = \min_{\alpha} \lambda_{\text{max}} (C P_\alpha C^\top)$  
          &
          $\displaystyle \| \mathcal{S}\|^2_{\infty, 2} = \lambda_{\text{max}} (C P C^\top)$ &
          $\displaystyle \| \mathcal{S}\|_\omega^2  = \min_{\tilde P} \lambda_{\text{max}} (C \tilde P C^\top) \; {\color{subsectioncolor} \star} $ \\
             &
          $\displaystyle \| \mathcal{S}\|_{\ast^\prime}^2  = \min_{\alpha} \lambda_{\text{max}} (B^\top Q_\alpha B) \; {\color{subsectioncolor} \star}$  
          &
          $\displaystyle \| \mathcal{S}\|^2_{2, i} = \lambda_{\text{max}} (B^\top Q B)$ &
          $\displaystyle \quad \quad \; \; \,  = \min_{\tilde Q} \lambda_{\text{max}} (B^\top \tilde Q B) \; {\color{subsectioncolor} \star}$ \\
          \midrule
          Trace norms &
          $\displaystyle \| \mathcal{S}\|_\varepsilon^2  = \min_{\alpha} \operatorname{trace} (C P_\alpha C^\top) \; {\color{subsectioncolor} \star}$  
          &
          $\displaystyle  \| \mathcal{S}\|^2_{\mathcal{H}_2} = \operatorname{trace} (C P C^\top)$ &
          $\displaystyle \| \mathcal{S}\|_\circ^2  = \min_{\tilde P} \operatorname{trace} (C \tilde P C^\top) \; {\color{subsectioncolor} \star}$ \\
             &
          $\displaystyle  \quad \quad \; \; = \min_{\alpha} \operatorname{trace} (B^\top Q_\alpha B) \; {\color{subsectioncolor} \star}$  
          &
          $\displaystyle \quad \quad \quad \; \:  = \operatorname{trace} (B^\top Q B)$ &
          $\displaystyle \| \mathcal{S}\|_{\circ^\prime}^2 = \min_{\tilde Q} \operatorname{trace} (B^\top \tilde Q B) \; {\color{subsectioncolor} \star}$ \\
          \midrule
          Gain estimates 
          & Peak-to-peak:
          & Energy-to-peak:
          & Integral-to-peak:
          \\[0.15cm] 
          &
          $\| \mathcal{S}\|_{\infty, \infty} \le  \| \mathcal{S}\|_\ast \le \| \mathcal{S}\|_\varepsilon$
          & $\| \mathcal{S}\|_{\infty, 2} \le \| \mathcal{S}\|_{\mathcal{H}_2}$ 
          & $\| \mathcal{S}\|_{\infty, 1} \le \| \mathcal{S}\|_\omega \le \| \mathcal{S}\|_\circ \; {\color{subsectioncolor} \star}$ \\[0.25cm] 
          & Impulse-to-integral:
          & Impulse-to-energy:
          & Impulse-to-peak:
          \\[0.15cm] 
          &
          $\| \mathcal{S}\|_{1, i} \le  \| \mathcal{S}\|_{\ast^\prime} \le \| \mathcal{S}\|_\varepsilon \; {\color{subsectioncolor} \star}$
          & $\| \mathcal{S}\|_{2, i} \le \| \mathcal{S} \| _{\mathcal{H}_2}$ 
          & $\| \mathcal{S}\|_{\infty, i} \le  \| \mathcal{S}\|_\omega \le \| \mathcal{S}\|_{\circ^\prime} \; {\color{subsectioncolor} \star}$ \\[0.15cm] 
          \bottomrule
          \end{tabular}
\end{table*}

It is well known that if $P, Q \succ 0$ are the solutions of \eqref{P2Q2}, i.e. controllability and observability Gramians of $\mathcal{S}$, then
\begin{equation*}
\operatorname{trace} (CPC^\top) = \operatorname{trace} (B^\top Q B).
\end{equation*}
We prove the similar properties of $P_\alpha, Q_\alpha$ and $\tilde P, \tilde Q$.

\begin{remark}
\label{remark1}
    Consider a system
    \begin{equation}
    \label{true_system}
            \dot x = \left(A+\frac{\alpha}{2}I\right)x + \frac{1}{\sqrt{\alpha}}Bu, \quad
             y = \frac{1}{\sqrt{\alpha}}Cx.
    \end{equation} 
    Notice that if $P_\alpha, Q_\alpha$ are the solutions of \eqref{Pinf}, \eqref{Q1}, then they are the controllability and observability Gramians of \eqref{true_system}. 
\end{remark}

\begin{proposition}
\label{pr:trace}
    Let $r<0$ be the value of the largest real part among all eigenvalues of $A$. If $\alpha \in (0, -2 r)$, then both \eqref{Pinf}, \eqref{Q1} admit positive definite solutions $P_\alpha, Q_\alpha \succ 0$, and 
    \begin{equation*}
         \operatorname{trace} (CP_\alpha C^\top) =  \operatorname{trace} (B^\top Q_\alpha B).
    \end{equation*}
    \vspace{-\baselineskip}
\end{proposition}

\begin{proposition}
\label{pr:max}
    If $\tilde P, \tilde Q \succ 0$ are subject to \eqref{Qinf}, \eqref{P1}, then 
    \begin{equation*}
        \min_{\tilde P} \lambda_{\text{max}} (C \tilde P C^\top) = \min_{\tilde Q} \lambda_{\text{max}} (B^\top \tilde Q B).
    \end{equation*}
    \vspace{-0.9\baselineskip}
\end{proposition}

    Note that in general
    \begin{equation*}
    \begin{aligned}
        \lambda_{\text{max}} (C P C^\top) &\ne  \lambda_{\text{max}} (B^\top  Q B), \\
        \lambda_{\text{max}} (C P_\alpha C^\top) & \ne \lambda_{\text{max}} (B^\top Q_\alpha B), \\
         \min_{\tilde P} \operatorname{trace} (C \tilde P C^\top) & \ne \min_{\tilde Q} \operatorname{trace} (B^\top \tilde Q B),
        \end{aligned}
    \end{equation*}
but all these equalities hold in SISO case, when $\operatorname{trace} = \lambda_{\text{max}}$.

\begin{remark}
In this section we have outlined the {{}duality relations} between $\mathcal{R}_p$ and $\mathcal{O}_q$ for $(p,q) = (1, \infty)$, $(2,2)$, $(\infty, 1)$. Note that all these pairs are Hölder conjugates.
\end{remark}

\subsection{System Norms}

Let $u \in \mathcal{F}^m(T)$, $y \in \mathcal{F}^k(T)$, $x(0)=0$. Define
\begin{equation*}
\| \mathcal{S}\|_{\infty, p} \vcentcolon = \sup_{T \ge 0} \: \max_{\|u\|_p \le 1} \, \|y\|_\infty .
\end{equation*}
For $p=2$ and $p=\infty$ this value is usually called ``energy-to-peak gain'' and ``peak-to-peak gain'' respectively (see \cite{b17-h2_norm}). We call it ``integral-to-peak  gain'' in case $p=1$.

Let $u(t) = u_0 \delta(t)$, $u_0 \in \mathbb{R}^m$, $y \in \mathcal{F}^k(T)$, $x(0)=0$, where $\delta (t)$ stands for the Dirac delta function. Define 
\begin{equation*}
\| \mathcal{S}\|_{q, i} \vcentcolon = \sup_{T \ge 0} \: \max_{|u_0| \le 1} \, \|y\|_q , 
\end{equation*}
where $i$ is a symbol that stands for ``impulse''. For $q=2$ and $q=\infty$ this value is usually called ``impulse-to-energy gain'' and ``impulse-to-peak gain'' respectively. We call it ``impulse-to-integral  gain'' in case $q=1$. 

It is known that if $P, Q \succ 0$ are the controllability and observability Gramians, i.e. the solutions of \eqref{P2Q2}, then
\begin{equation*}
  \| \mathcal{S}\|^2_{\infty, 2} = \lambda_{\text{max}} (C P C^\top), \quad \| \mathcal{S}\|^2_{2, i} =  \lambda_{\text{max}} (B^\top  Q B), 
\end{equation*}
\vspace{-\baselineskip}
\begin{equation*}
  \| \mathcal{S}\|^2_{\mathcal{H}_2} = \operatorname{trace} (CPC^\top) = \operatorname{trace} (B^\top Q B),
\end{equation*}
where $\| \mathcal{S}\|_{\mathcal{H}_2}$ is the usual $\mathcal{H}_2$-norm of $\mathcal{S}$ (see \cite{b17-h2_norm}).

Let $P_\alpha, Q_\alpha \succ 0$ be the solutions of \eqref{Pinf}, \eqref{Q1}. The $\ast$-norm, which was studied in \cite{b5-abedor}, \cite{b6-starcapture}, is typically defined as 
\begin{equation*}
  \| \mathcal{S}\|_\ast^2 \vcentcolon = \min_{\alpha} \lambda_{\text{max}} (C P_\alpha C^\top).
\end{equation*}
Define its dual counterpart, the $\ast^\prime$-norm, as
\begin{equation*}
  \| \mathcal{S}\|_{\ast^\prime}^2 \vcentcolon = \min_{\alpha} \lambda_{\text{max}} (B^\top Q_\alpha B).
\end{equation*}
We introduce the family of $\varepsilon(\alpha)$-norms and the $\varepsilon$-norm, defined as
\begin{equation*}
\| \mathcal{S}\|_{\varepsilon (\alpha)}^2 \vcentcolon =  \operatorname{trace} (CP_\alpha C^\top) =  \operatorname{trace} (B^\top Q_\alpha B), 
\end{equation*}
\begin{equation*}
  \| \mathcal{S}\|_\varepsilon \vcentcolon = \min_\alpha \| \mathcal{S}\|_{\varepsilon (\alpha)}.
\end{equation*}
The $\varepsilon(\alpha)$-norm is defined for $\alpha \in (0, -2r)$ by Proposition \ref{pr:trace}. The $\varepsilon$-norm is well-defined (the minimum is achieved) by the convexity of $\varphi : \alpha \mapsto \operatorname{trace} CP_\alpha C^\top$ proven in\cite{b1-nazin}.


Let $\tilde P, \tilde Q \succ 0$ be the solutions of \eqref{Qinf}, \eqref{P1}. Define
\begin{equation*}
  \| \mathcal{S}\|_\omega^2 \vcentcolon = \min_{\tilde P} \lambda_{\text{max}} (C \tilde P C^\top) = \min_{\tilde Q} \lambda_{\text{max}} (B^\top \tilde Q B),
\end{equation*}
\vspace{-\baselineskip}
\begin{equation*}
  \| \mathcal{S}\|_\circ^2  \vcentcolon =  \min_{\tilde P} \operatorname{trace} (C \tilde P C^\top), \; 
  \| \mathcal{S}\|_{\circ^\prime}^2  \vcentcolon =  \min_{\tilde Q} \operatorname{trace} (B^\top \tilde Q B),
\end{equation*}
where the $\omega$-norm is well-defined by Proposition \ref{pr:max}.

It is natural to compare newly introduced norms of ``largest eigenvalue'' and ``trace'' type with the system gains $\| \mathcal{S}\|_{\infty,p}$ and $\| \mathcal{S}\|_{q,i}$. We do it by means of Propositions \ref{pr:Pinequality}, \ref{pr:Qinequality} and Theorem \ref{thm:norms}.

\begin{proposition}
\label{pr:Pinequality}
    If $\mathcal{P} \succ 0$ is a matrix of an outer ellipsoidal approximation for $\mathcal{R}_p$, i.e.
    \begin{equation*}
        \mathcal{R}_p \subset  \left\{ x \mid x^\top \mathcal{P}^{-1} x \le 1 \right\} ,
    \end{equation*}
then 
\begin{equation*}
    \|\mathcal{S}\|_{\infty,p}^2 \le \lambda_{\text{max}} (C \mathcal{P} C^\top) \le \operatorname{trace} (C \mathcal{P} C^\top).
    \end{equation*}
    \vspace{-\baselineskip}
\end{proposition}

Note that the geometrical meaning of $\lambda_{\text{max}} (C \mathcal{P} C^\top)$ is the square of the largest semiaxis of an ellipsoid
\begin{equation*}
    \left\{ y \mid y^\top (C\mathcal{P} C^\top)^{-1} y \le 1 \right\},
\end{equation*}
while $\operatorname{trace} (C \mathcal{P} C^\top)$ is the sum of the squares of its semiaxes. Hence, when we aim at minimizing these values, we are in fact trying to make an outer ellipsoidal approximation of a {{}easily reachable} set $\mathcal{R}_p$ as small as possible, at the same time tightening the estimate of  $\|\mathcal{S}\|_{\infty,p}$.

\begin{proposition}
\label{pr:Qinequality}
    If $\mathcal{Q} \succ 0$ is a matrix of an inner ellipsoidal approximation for $\mathcal{O}_q$, i.e.
    \begin{equation*}
        \left\{ x \mid x^\top \mathcal{Q} x \le 1 \right\} \subset \mathcal{O}_q    ,
    \end{equation*}
then 
\begin{equation*}
    \|\mathcal{S}\|_{q,i}^2 \le \lambda_{\text{max}} (B^\top  \mathcal{Q} B) \le \operatorname{trace} (B^\top  \mathcal{Q} B).
    \end{equation*}
    \vspace{-\baselineskip}
\end{proposition}

Note that the geometrical meaning of $\lambda_{\text{max}} (B^\top  \mathcal{Q} B)$ is the inverse square of the smallest semiaxis of an ellipsoid
\begin{equation*}
    \left\{ u \mid u^\top (B^\top  \mathcal{Q} B) u \le 1 \right\},
\end{equation*}
while $\operatorname{trace} (B^\top  \mathcal{Q} B)$ is the sum of the inverse squares of its semiaxes. Hence, when we aim at minimizing these values, we are in fact trying to make an inner ellipsoidal approximation of a {{}hardly observable} set $\mathcal{O}_q$ as large as possible, at the same time tightening the estimate of  $\|\mathcal{S}\|_{q,i}$. 

Now we can establish the theorem, that provides the estimates for the system gains in terms of the system norms with ellipsoidal geometrical meaning. 

\begin{theorem}
\label{thm:norms}
The following estimates hold
\begin{alignat*}{7}
&   \text{\textbullet \; Peak-to-peak gain:} && \| \mathcal{S}\|_{\infty, \infty} && \le \; && \| \mathcal{S}\|_\ast && \le \; && \| \mathcal{S}\|_\varepsilon && \, ;\\
& \text{\textbullet \; Impulse-to-integral gain:} \quad \; && \| \mathcal{S}\|_{1, i} && \le  && \| \mathcal{S}\|_{\ast^\prime} && \le && \| \mathcal{S}\|_\varepsilon && \, ;\\
& \text{\textbullet \; Integral-to-peak gain:} && \| \mathcal{S}\|_{\infty, 1} && \le && \| \mathcal{S}\|_\omega && \le && \| \mathcal{S}\|_\circ && \, ; \\
& \text{\textbullet \; Impulse-to-peak gain:} && \| \mathcal{S}\|_{\infty, i} && \le  && \| \mathcal{S}\|_\omega && \le && \|  \mathcal{S}\|_{\circ^\prime} && \, .
\end{alignat*}
\vspace{-\baselineskip}
\end{theorem}

It is important to note that the $\varepsilon$-norm and the $\omega$-norm are of significant value as they are related to both external approximations of {{}easily reachable} sets and internal approximations of {{}hardly observable} sets. While the $\varepsilon$-norm may not be the sharpest estimate of the peak-to-peak gain, it provides a natural and symmetric estimate of both the peak-to-peak and impulse-to-integral gains. The $\omega$-norm also benefits from its symmetry, in addition to being a less conservative estimate compared to $\circ$ and $\circ^\prime$ norms. {{} Moreover, the $\varepsilon$-norm is directly related to the magnitudes of all the semiaxes of the corresponding approximating ellipsoid. In contrast, the $\ast$-norm only captures the length of its largest semiaxis. Therefore, if the objective is to minimize the system's responsiveness across all directions in the output space (and not only the worst-case one), the minimization of the $\varepsilon$-norm is more advantageous. For this reason, the upcoming sections on synthesis will focus on minimizing the $\varepsilon$-norm of the closed-loop system.}

\begin{figure*}[t!]
\begin{subfigure}{0.5\textwidth}
\centering
\includegraphics[width = 0.92\textwidth]{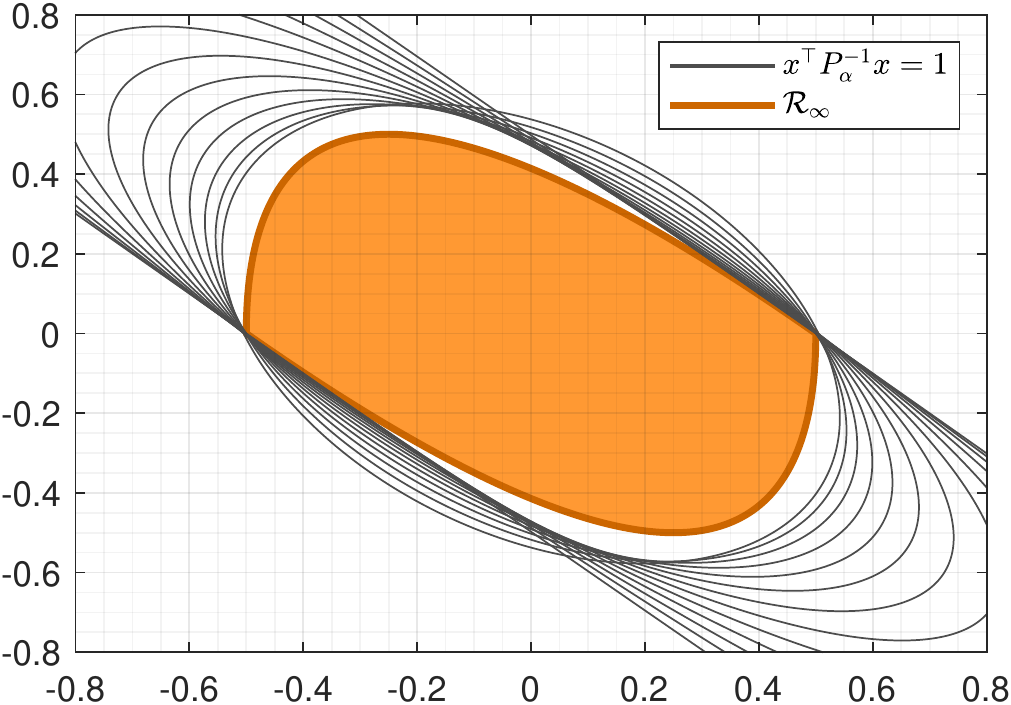}
\caption{}
\vspace*{2mm}
\end{subfigure}
\hfill
\begin{subfigure}{0.5\textwidth}
\centering
\includegraphics[width = 0.92\textwidth]{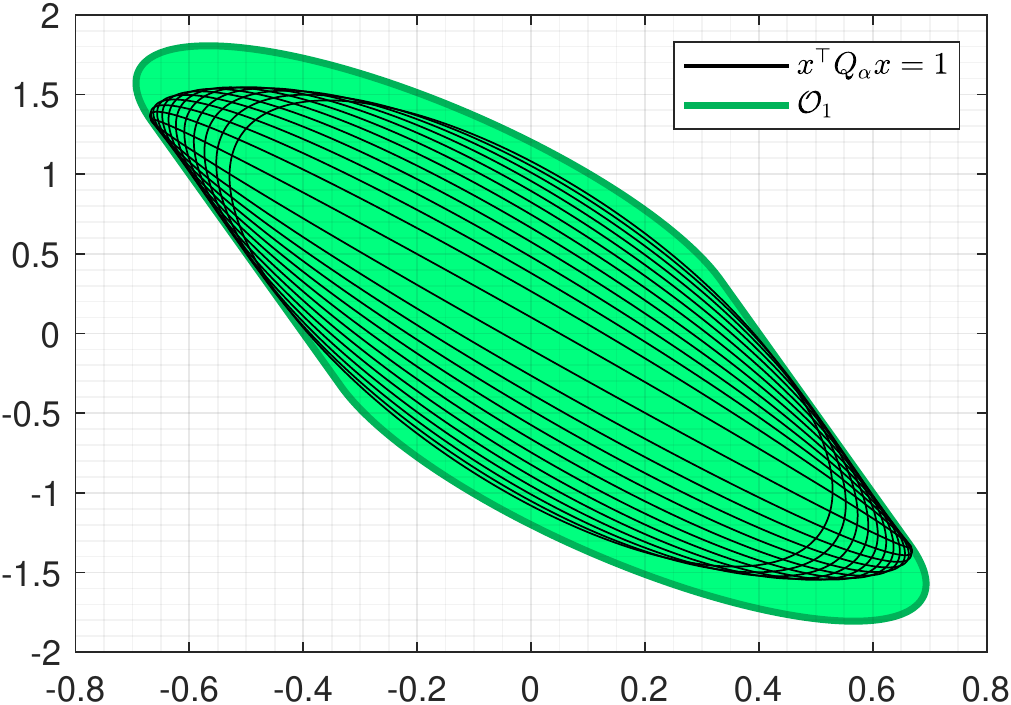}
\caption{}
\vspace*{2mm}
\end{subfigure}
\bigskip
\begin{subfigure}{0.5\textwidth}
\centering
\includegraphics[width = 0.92\textwidth]{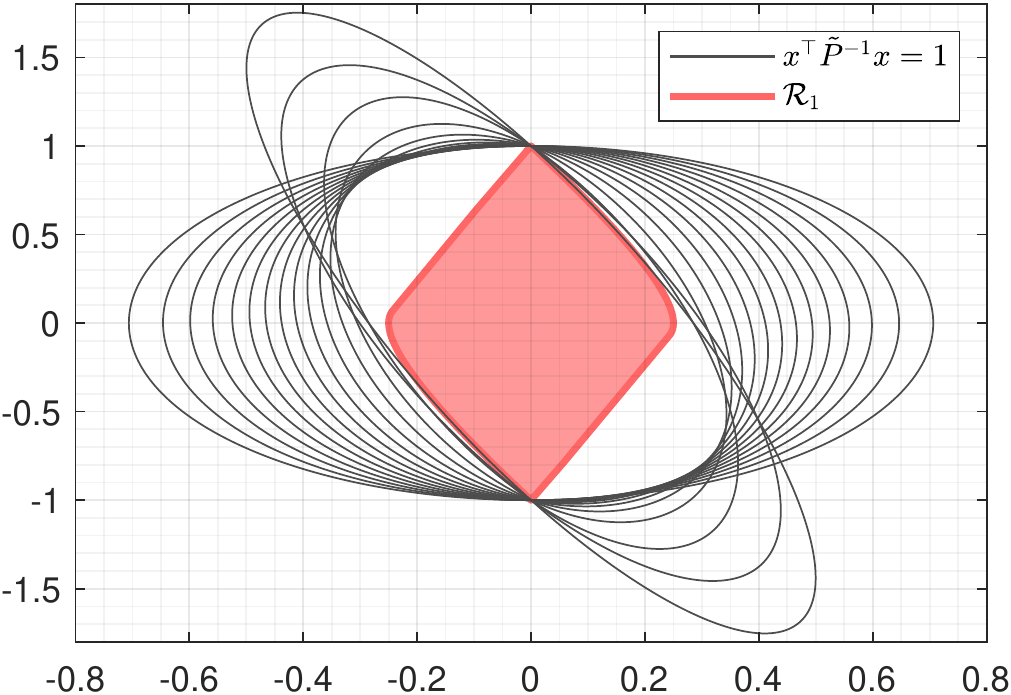}
\caption{}
\end{subfigure}
\hfill
\begin{subfigure}{0.5\textwidth}
\centering
\includegraphics[width = 0.92\textwidth]{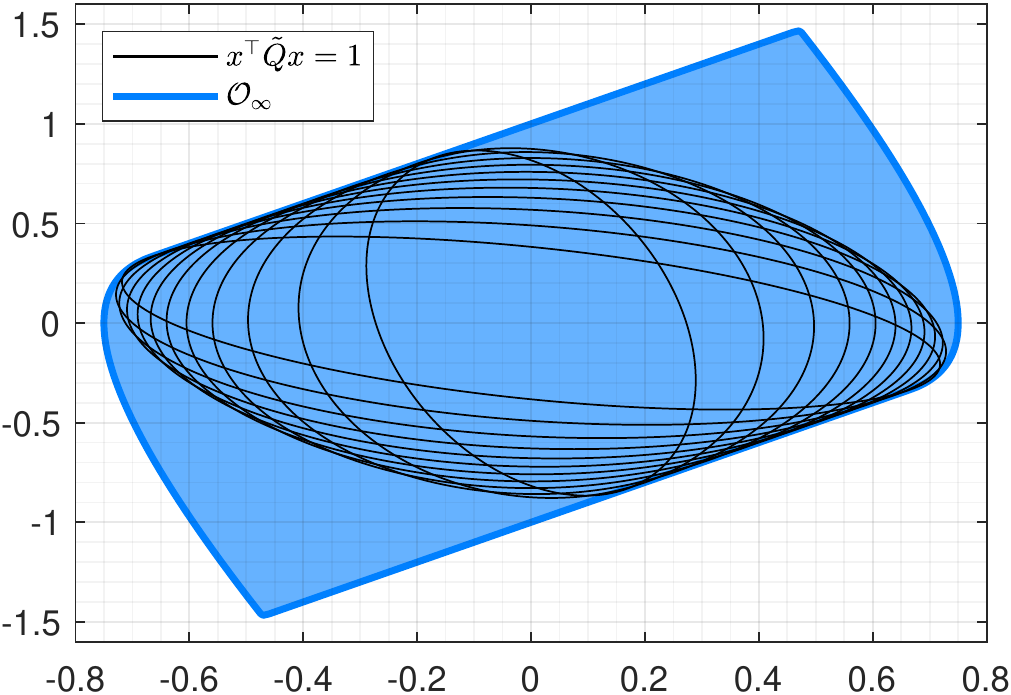}
\caption{}
\end{subfigure}
\vspace*{-8mm}
\caption{Illustrations for the example from section \ref{example}. Filled areas indicate {{}easily reachable} sets with respect to (a) $\infty$-norm, (c) $1$-norm, and {{}hardly observable} sets with respect to (b) $1$-norm, (d) $\infty$-norm. Black curves correspond to their ellipsoidal approximations obtained from \eqref{Pinf}-\eqref{P1}. } 
\label{fig1}
\end{figure*}

\subsection{Illustrative example} \label{example}

Consider a system $\mathcal{S}$ with matrices
\begin{equation*}
    A = \begin{bmatrix}
    0 & 1 \\ -2 & -3
    \end{bmatrix}, \quad 
    B = \begin{bmatrix}
        0 \\ 1
    \end{bmatrix}, \quad
    C = \begin{bmatrix}
        1 & -1
    \end{bmatrix}.
\end{equation*}

Fig. \ref{fig1} shows the sets $\mathcal{R}_1$, $\mathcal{R}_\infty$,  $\mathcal{O}_1$, $\mathcal{O}_\infty$ for this system, along with their ellipsoidal approximations. The sets $\mathcal{R}_1$ and $\mathcal{R}_\infty$ were constructed using the support function method (refer to  \cite{b7-support1}, \cite{b8-support2} for further details), while the sets $\mathcal{O}_1$ and $\mathcal{O}_\infty$ were determined through direct calculations. Ellipsoidal approximations were obtained from \eqref{Pinf}-\eqref{P1}. In accordance with Theorems \ref{thm:Pinf}-\ref{thm:P1}, approximations of {{}easily reachable} sets are external, while {{}hardly observable} sets are approximated from the inside. System gains and their estimates are
\begin{alignat*}{5}
    &\| \mathcal S \|_{\infty,\infty} && = 0.833 \; &&   \le \; 0.914 = \| \mathcal S\|_{\ast} && = \| \mathcal S\|_\varepsilon ,   \\   
    &\| \mathcal S \|_{1,i} && = 0.833 &&  \le \; 0.914 = \| \mathcal S\|_{\ast^\prime} && = \| \mathcal S\|_\varepsilon , \\
     &\| \mathcal S \|_{\infty,1} && = 1 &&  \le \;  1.144 =  \| \mathcal S\|_{\omega}  && = \| \mathcal S\|_\circ  ,\\ 
    &\| \mathcal S \|_{\infty,i} && = 1 &&  \le \; 1.144 =  \| \mathcal S\|_{\omega}    && = \| \mathcal S\|_{\circ^\prime} ,
\end{alignat*}
where the equality between some norms is due to the fact that this system is SISO. The value of $\| \mathcal{S}\|_\varepsilon = \left( \min CP_\alpha C^\top \right)^{1/2}$ was obtained at $\alpha = 0.67$. Note that Theorem \ref{thm:norms} holds.

\section{Synthesis: State-Feedback and Filtering}
\label{S3}

In this section we propose a way to design the optimal state-feedback controller, as well as the optimal observer (filter) with respect to the $\varepsilon$-norm. 

A natural way to understand the optimality of the proposed solutions is to consider the situation of bounded external disturbances, when one wants to minimize the peak-to-peak gain $\|\mathcal{S}\|_{\infty, \infty}$ of a closed-loop system, or tries to achieve the smallest size of set $\mathcal{R}_\infty$ of states reachable by bounded disturbances. However, as the $\varepsilon$-norm is also an upper bound for the impulse-to-integral gain $\|\mathcal{S}\|_{1, i}$, the solutions are implicitly optimal with respect to this dual criterion as well.

\subsection{Optimal State-Feedback with respect to \texorpdfstring{$\varepsilon$}{ε}-norm}

Consider a linear time-invariant plant
\begin{equation}
\label{plant1}
    \left\{ 
    \begin{aligned}
        &\dot x = Ax + B u + B_w w, \\
        &z = C x + D u,
    \end{aligned}
    \right.
\end{equation}
where {{} $x(t) \in \mathbb{R}^n$, $u(t) \in \mathbb{R}^m$, $w(t) \in \mathbb{R}^{\bar m}$, $z(t) \in \mathbb{R}^k$, and $A$, $B$, $B_w$, $C$, $D$ are real matrices of corresponding sizes.} Here we regard $w$ as the external disturbance signal and $z$ as the regulated output. We make the standard assumptions that $(A,B)$ is stabilizable, $(C,A)$ is observable, $C^\top D = 0$, and $D^\top D$ is invertible. 

Consider a linear static feedback controller of the form
\begin{equation}
\label{ctrl1} 
    u = Kx,
\end{equation}
where $K \in \mathbb{R}^{m \times n}$. Denote the closed-loop system \eqref{plant1}-\eqref{ctrl1} as $\mathcal{S}_K$. Regard $\mathcal{S}_K$ as the system with the input $w$ and the output $z$. {{}Consider the following optimal control problem. 

\begin{problem} \label{problem1}
     Find the optimal controller matrix $K$ that minimizes the  $\varepsilon$-norm of the closed-loop system: $\| \mathcal{S}_{K} \|_\varepsilon \to \min.$
\end{problem}
}

\begin{proposition}
\label{pr:Q-solvable}
    If $(A,B)$ is stabilizable, $(C,A)$ is observable, then for each $\alpha > 0$ the equation
    \begin{equation} \begin{multlined}[c][.85\displaywidth]
        Q_\alpha A + A^\top Q_\alpha  + \alpha Q_\alpha  \\ - \alpha Q_\alpha B(D^\top D)^{-1} B^\top Q_\alpha + \frac{1}{\alpha} C^\top C = 0
        \label{ricQ}
    \end{multlined} \end{equation} 
    admits the unique positive definite solution $Q_\alpha \succ 0$.
\end{proposition}

\begin{theorem}
\label{thm:K}
    Let $Q_\alpha$ be the positive definite solution of \eqref{ricQ}. The controller \eqref{ctrl1} with matrix
    \begin{equation}
    \label{ricQ-K}
        K = - \alpha (D^\top D)^{-1}B^\top Q_{\alpha}
    \end{equation}
    renders the system $\mathcal{S}_{K \!}$ stable and guarantees that its $\varepsilon(\alpha)$-norm possesses the smallest possible value, which is equal to
    \begin{equation*}
        \| \mathcal{S}_K \|_{\varepsilon(\alpha)} = \sqrt{\operatorname{trace} (B_w^{\top} Q_{\alpha} B_w^{\phantom \top})}. 
    \end{equation*}
     \vspace{-\baselineskip}
\end{theorem}

According to Theorem \ref{thm:K}, finding the {{}solution to Problem~\ref{problem1}} involves iterating the parameter $\alpha \in (0, \infty)$ and selecting the one that minimizes $\operatorname{trace} (B_w^{\top} Q_{\alpha} B_w^{\phantom \top})$. The corresponding controller gain \eqref{ricQ-K} will be optimal in terms of the $\varepsilon$-norm.

Previous studies \cite{b1-nazin}, \cite{b3-poznyak} solved the task $\| \mathcal{S}_K \|_\varepsilon \to \min$ as an optimization problem  for each fixed $\alpha$. The solution had to be obtained by iterating the parameter and solving the system of LMIs with $\frac{1}{2}n(n+1)+\frac{1}{2}m(m+1)+mn$ variables on every iteration. In contrast to that, the proposed solution is in the form of the parameter depending Riccati equation \eqref{ricQ} with only $\frac{1}{2}n(n+1)$ variables. It requires iterating the parameter $\alpha$ as well, but instead of a system of LMIs, one has to solve a specific Riccati equation on every iteration. The advantages of this approach is threefold: 

\begin{itemize}
    \item With standard tools (\texttt{cvx} software and MATLAB function \texttt{are}), Riccati equation \eqref{ricQ} can be solved faster than the system of LMIs from \cite{b1-nazin}, \cite{b3-poznyak}. It becomes particularly evident when iterating the parameter $\alpha$. 
    
    \item Theorem \ref{thm:K} will be used in the sequel to solve the optimal output-feedback control problem in Section \ref{S4}. 
    
    \item It links the obtained results with the famous $\mathcal{H}_2$-control. However, there are major differences between the two (see Section \ref{S5B}).
\end{itemize}

 The following proposition guarantees that the best $\alpha$ (i.e. the one that leads to the smallest $\varepsilon$-norm) is always achieved away form $0$.

\begin{proposition}
\label{pr:notzero}
    If $Q_\alpha \succ 0$ is given by \eqref{ricQ}, and
    \begin{equation*}
        \hat \alpha \vcentcolon = \underset{\alpha \in (0, \infty)} {\operatorname{arginf}} \operatorname{trace} (B_w^{\top} Q_{\alpha} B_w^{\phantom \top}),
    \end{equation*}
    then $\hat \alpha \ne 0$ (i.e. either $\hat \alpha \in (0, \infty)$, or $\hat \alpha = \infty$). 
\end{proposition}

\begin{remark}
    In \cite{b13-khlebnikov_filter} it was conjectured that $\alpha \mapsto \| \mathcal{S}_K \|_{\varepsilon(\alpha)}^2$ is always a convex function, which would make the search for $\hat \alpha$ a lot easier. We give a simple counterexample to this conjecture. Consider the plant \eqref{plant1} with matrices
    \begin{equation*}
    A = \begin{bmatrix}
          0 & 1 & 0 \\ 0 & 0 & 1 \\ 1 & 0 & 1
    \end{bmatrix}\!,
    \;
    B = \begin{bmatrix}
        0 \\ 1 \\ 1
    \end{bmatrix}\!,
    \;
    B_w = \begin{bmatrix}
        2 \\ 1 \\ 0
    \end{bmatrix}\!, 
    \;
    C^\top = \begin{bmatrix}
          1 & 0 \\ 0 & 0 \\ 10 & 0
    \end{bmatrix}\!,
    \end{equation*}
    \vspace{-0.5\baselineskip}
\end{remark}
and $D = \begin{bmatrix} 0 & 1 \end{bmatrix}^\top$. By numerical study one can show that the function under discussion has at least two local minima, namely $\alpha \approx 0.09$ and $\alpha \approx 2.06$. 

\begin{remark}
\label{remark4}
Determining necessary and sufficient conditions for $\hat \alpha$ to be finite is of interest. For $m = \bar m = 1$, linear independence of $B$ and $B_w$ seems sufficient, but we leave the exact formulation as an open problem for future research.
\end{remark}

\subsection{Optimal Filtering with respect to \texorpdfstring{$\varepsilon$}{ε}-norm}

Consider a linear time-invariant plant
\begin{equation}
\label{plant2}
    \left\{ 
    \begin{aligned}
        &\dot x = Ax + B w, \\
        &y = C x + D w, \quad \quad
    \end{aligned}
    \right.
\end{equation}
where {{} $x(t) \in \mathbb{R}^n$, $w(t) \in \mathbb{R}^{m}$, $y(t) \in \mathbb{R}^k$, and $A$, $B$, $C$, $D$ are real matrices of corresponding sizes.} Here we regard $w$ as the external disturbance signal and $y$ as the measured output. We make the standard assumptions that $(C,A)$ is detectable, $(A,B)$ is controllable, $B D^\top = 0$ and $D D^\top$ is invertible. 

Consider a linear time-invariant observer of the form
\begin{equation}
\label{ctrl2}
    \left\{ 
    \begin{aligned}
        &\dot {\hat x} = A\hat x + L(\hat y - y), \\
        & \hat y = C \hat x,
    \end{aligned}
    \right.
\end{equation}
and the observer error 
\begin{equation}
\label{ctrl3}
    z = C_z (x-\hat x), \quad \; \;
\end{equation}
where $L \in \mathbb{R}^{n \times k}$, $C_z \in \mathbb{R}^{\bar{k} \times n}$. Denote the closed-loop system \eqref{plant2}-\eqref{ctrl3} as $\mathcal{S}_L$. Regard $\mathcal{S}_L$ as the system with the input $w$ and the output $z$. {{}Consider the following optimal observer design problem. 

\begin{problem} \label{problem2}
     Find the optimal observer gain $L$ that minimizes the  $\varepsilon$-norm of the closed-loop system: $\| \mathcal{S}_{L} \|_\varepsilon \to \min.$
\end{problem}
}

\begin{proposition}
\label{pr:P-solvable}
    If $(A,B)$ is controllable, $(C,A)$ is detectable, then for each $\alpha > 0$ the equation
    \begin{equation} \begin{multlined}[c][.85\displaywidth]
         A P_\alpha + P_\alpha A^\top  + \alpha P_\alpha  \\ - \alpha P_\alpha C^\top (D D^\top)^{-1} C P_\alpha + \frac{1}{\alpha} B B^\top = 0
        \label{ricP}
    \end{multlined} \end{equation} 
    admits the unique positive definite solution $P_\alpha \succ 0$.
\end{proposition}

\begin{theorem}
\label{thm:L}
    Let $P_\alpha$ be the positive definite solution of \eqref{ricP}.
  The observer \eqref{ctrl2} with matrix
    \begin{equation}
    \label{ricP-L}
        L = - \alpha P_{\alpha} C^\top (D D^\top)^{-1} 
    \end{equation}
    renders the system $\mathcal{S}_L$ stable and guarantees that its $\varepsilon (\alpha)$-norm possesses the smallest possible value, which is equal to
    \begin{equation*}
        \| \mathcal{S}_L \|_{\varepsilon(\alpha)} = \sqrt{\operatorname{trace} (C_z^{\phantom{;}} P_{ \alpha} C_z^\top) }. 
    \end{equation*}
     \vspace{-\baselineskip}
\end{theorem}

According to Theorem \ref{thm:L}, finding the {{}solution to Problem~\ref{problem2}} involves iterating the parameter $\alpha \in (0, \infty)$ and selecting the one that minimizes $\operatorname{trace} (C_z^{\phantom{;}} P_{ \alpha} C_z^\top)$. The corresponding observer gain \eqref{ricP-L} will be optimal in terms of the $\varepsilon$-norm.

In previous work \cite{b13-khlebnikov_filter} the task $\| \mathcal{S}_L \|_\varepsilon \to \min$ was solved as an optimization problem with the constraints given by LMIs for every fixed value of the parameter $\alpha$. The solution had to be obtained by iterating the parameter and solving the system of LMIs with $n(n+1)+nk$ variables on every iteration. In contrast to that, the proposed solution is in the form of the parameter depending Riccati equation \eqref{ricP} with only $\frac{1}{2}n(n+1)$ variables. It requires iterating the parameter $\alpha$ as well, but instead of a system of LMIs, one has to solve a specific Riccati equation on every iteration. This solution has the same advantages as the state-feedback one, and will {{}also} be used in the sequel.

\begin{proposition}
\label{pr:notzero2}
    If $P_\alpha \succ 0$ is given by \eqref{ricP}, and
    \begin{equation*}
        \hat \alpha \vcentcolon = \underset{\alpha \in (0,\infty)} {\operatorname{arginf}} \operatorname{trace} (C_z^{\phantom{;}} P_{ \alpha} C_z^\top),
    \end{equation*}
    then $\hat \alpha \ne 0$ (i.e. either $\hat \alpha \in (0, \infty)$, or $\hat \alpha = \infty$). 
\end{proposition}

\begin{remark}

Determining necessary and sufficient conditions for $\hat \alpha$ to be finite is of interest. For $k = \bar k = 1$, linear independence of $C$ and $C_z$ seems sufficient, but we leave the exact formulation as an open problem for future research.
\end{remark}

\section{Synthesis: Output-Feedback Control}
\label{S4}

Consider a linear time-invariant plant
\begin{equation}
\label{plant_main}
   \left\{ 
    \begin{aligned}
        & \dot x = Ax + B_1 w + B_2 u, \\
        & y = C_1 x + D_1 w, \\
        & z = C_2 x + D_2 u,
    \end{aligned}
    \right. \; \; \; 
\end{equation}
where {{} $x(t) \in \mathbb{R}^n$, $u(t) \in \mathbb{R}^m$, $y(t) \in \mathbb{R}^k$, $w(t) \in \mathbb{R}^{\bar m}$, $z(t) \in \mathbb{R}^{\bar k}$, and $A$, $B_i$, $C_i$, $D_i$ are real matrices of corresponding sizes.} Regard $u$ as the control input, $w$ as the external disturbance, $y$ as the measured output, and $z$ as the regulated output. Assume that $(A,B_2)$ is stabilizable, $(C_1, A)$ is detectable, $(A,B_1)$ is controllable, $(C_2,A)$ is observable, $B_1 D_1^\top = 0$, $C_2^{\top} D_2 = 0$, and both  $D_1 D_1^\top$ and $D_2^\top D_2$ are invertible.

Consider a linear time-invariant controller of the form
\begin{equation}
\label{ctrl_main}
   \quad \left\{ 
    \begin{aligned}
        & \dot{\hat x} = A \hat{x} + B_2 u + L(\hat{y}-y), \\
        & \hat{y} = C_1 \hat x, \\
        & u = K \hat x,
    \end{aligned}
    \right.
\end{equation}
where $K \in \mathbb{R}^{m \times n}$, $L \in \mathbb{R}^{n \times k}$. Note that this controller structure  is the classic output-feedback controller, that combines equations \eqref{ctrl1} and \eqref{ctrl2}.

Denote the closed-loop system \eqref{plant_main}-\eqref{ctrl_main} as $\mathcal{S}_{KL}$. Regard $\mathcal{S}_{KL}$ as the system with the input $w$ and the output $z$. {{}Consider the following optimal control problem. 

\begin{problem} \label{problem3}
     Find the optimal controller and observer matrices $K$, $L$ that minimize the  $\varepsilon$-norm of the closed-loop system:
\begin{equation*}
    \| \mathcal{S}_{KL} \|_\varepsilon \to \min.
\end{equation*}
\vspace{-\baselineskip}
\end{problem}
}

\begin{theorem}[(Main result)]
    \label{thm:KL}
    Let $Q_\alpha$ be the positive definite solution of Riccati equation
\begin{equation*}
\begin{multlined}[c][.9\displaywidth]
        Q_\alpha A + A^\top Q_\alpha  + \alpha Q_\alpha  \\ - \alpha Q_\alpha B_2(D_2^\top D_2)^{-1} B_2^\top Q_\alpha + \frac{1}{\alpha} C_2^\top C_2 = 0,
\end{multlined}
\end{equation*}
and $P_\alpha$ be the positive definite solution of Riccati equation
\begin{equation*}
\begin{multlined}[c][.9\displaywidth]
        A P_\alpha + P_\alpha A^\top  + \alpha P_\alpha  \\ - \alpha P_\alpha C_1^\top (D_1 D_1^\top)^{-1} C_1 P_\alpha + \frac{1}{\alpha} B_1 B_1^\top = 0.
\end{multlined}
\end{equation*}
The controller \eqref{ctrl_main} with matrices 
\begin{equation*}
    K = - \alpha (D_2^\top D_2)^{-1}B_2^\top Q_\alpha, \quad
    L = - \alpha P_\alpha C_1^\top (D_1 D_1^\top)^{-1}
\end{equation*}
renders the system $\mathcal{S}_{KL}$ stable and guarantees that its $\varepsilon(\alpha)$-norm possesses the smallest possible value, which is 
\begin{equation*} 
\begin{aligned}
    \| S_{KL} \|_{\varepsilon(\alpha)} \! & = \! \left( \operatorname{trace} (B_1^\top  Q_{\alpha} B_1) + \operatorname{trace} (D_2 K P_{ \alpha} K^\top D_2^\top) \right)^{\scriptscriptstyle \! 1/2} \\
    & = \! \left( \operatorname{trace} (C_2  P_{\alpha} C_2^\top) + \operatorname{trace} (D_1^\top L^\top Q_{ \alpha} L D_1)\right)^{\scriptscriptstyle \! 1/2}.
\end{aligned}
\end{equation*}
     \vspace{0cm}
\end{theorem}

According to Theorem \ref{thm:KL}, finding the {{}solution to Problem~\ref{problem3}} involves iterating the parameter $\alpha \in (0, \infty)$ and selecting the one that minimizes $\| \mathcal{S}_{KL} \|_{\varepsilon (\alpha)}$. {{}Note that due to Propositions \ref{pr:Q-solvable} and \ref{pr:P-solvable} the existence of the corresponding $Q_\alpha$ and $P_\alpha$ is guaranteed for all $\alpha > 0$}.

\begin{remark}
\label{remark6}
    Numerical simulations suggest that the function $\alpha \mapsto \left( \operatorname{trace} (B_1^\top  Q_\alpha B_1) + \operatorname{trace} (D_2 K P_\alpha K^\top D_2^\top) \right) = \left( \operatorname{trace} (C_2  P_\alpha C_2^\top) + \operatorname{trace} (D_1^\top L^\top Q_\alpha L D_1) \right)$ is strictly convex and that the value
    \begin{equation*}
        \hat \alpha = \underset{\alpha > 0}{\operatorname{argmin}} \| \mathcal{S}_{KL} \|_{\varepsilon (\alpha)}
    \end{equation*}
    is always finite and nonzero. However, a rigorous justification of this fact remains an open problem. 
\end{remark}

\section{Comparison with Known Results}
\label{S5}

\subsection{Comparison with the invariant ellipsoids method}

Before this work, the only known solution for the problem $\| \mathcal{S}_{KL} \|_\varepsilon \to \min$ for \eqref{plant_main}, \eqref{ctrl_main} was presented in \cite{b2-topunov} {{}and then refined in \cite{b13-khlebnikov_filter}}. Though it was acknowledged that this solution only provided sub-optimal results, it became a generally accepted method and its extended version was discussed in the well-known monograph \cite{b3-poznyak}.

In the present paper the optimality of the solution proposed in Theorem \ref{thm:KL} is proved.  From a theoretical perspective, this means that there is no other possible choice of matrices $K$ and $L$ that can further minimize the value of $\| \mathcal{S}_{KL} \|_\varepsilon$. However, it is still of interest to compare this proposed optimal solution with the previously known sub-optimal one, to show the difference between the two.

Consider the plant \eqref{plant_main} with matrices $ A = \begin{bmatrix}
        0 & 1 \\ \beta & 0
    \end{bmatrix}$, $\beta \in \mathbb{R}$,
\begin{alignat*}{3}
    & B_1 = \begin{bmatrix}
        1 & 0 & 0 \\
        0 & 1 & 0
    \end{bmatrix}, \; \;
    && B_2 = \begin{bmatrix}
        0 \\ 1
    \end{bmatrix}, \; \;
    && D_1 = \begin{bmatrix}
        0 & 0 & 1
    \end{bmatrix},
    \\
    & C_2 = \begin{bmatrix}
        1 & 0 & 0 \\
        0 & 1 & 0
    \end{bmatrix}^\top, \; \;
    && C_1  = \begin{bmatrix}
        1 \\ 0
    \end{bmatrix}^\top, \; \;
    && D_2 = \begin{bmatrix}
        0 & 0 & 1
    \end{bmatrix}^\top \! \!.
\end{alignat*}
We deliberately selected a simple plant for this comparison to make it easily reproducible and to clearly demonstrate the differences between the two algorithms.

\begin{figure}[t]
    \centering
    \includegraphics[height = 6.3cm]{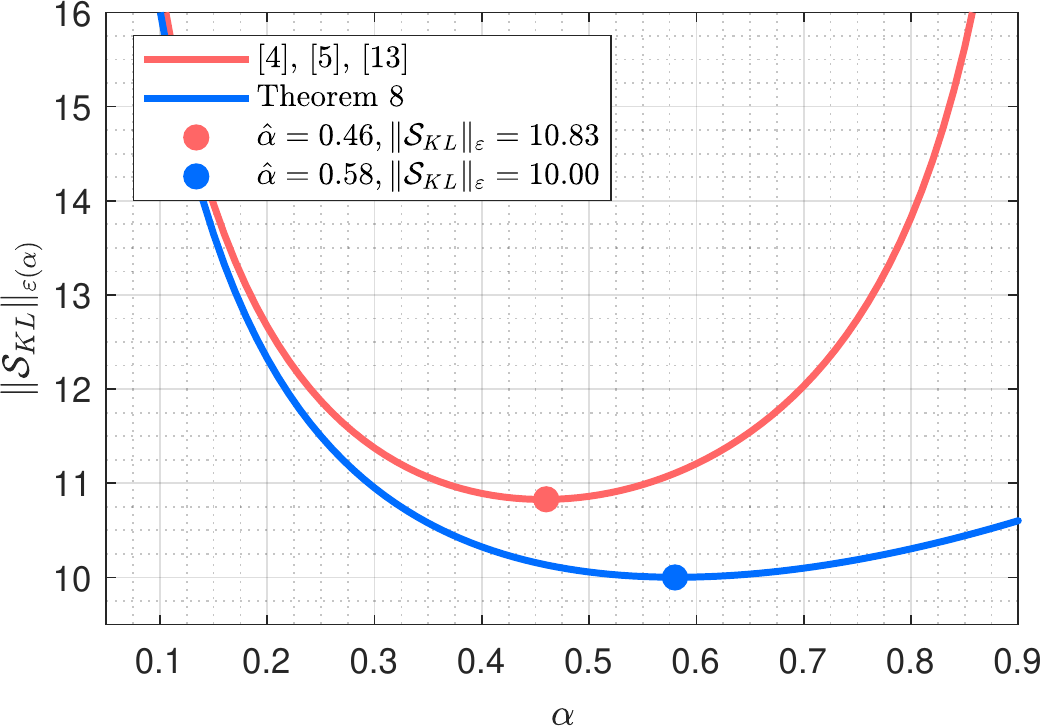}
    \caption{A comparison of the closed-loop system's $\varepsilon(\alpha)$-norm with {{}$\beta = 0.3$ between the sub-optimal controller \cite{b13-khlebnikov_filter}, \cite{b2-topunov}, \cite{b3-poznyak}} and the optimal controller proposed in Theorem \ref{thm:KL}.}
    \label{fig2}
\end{figure}

\begin{figure}[t]
    \centering
    \includegraphics[height = 6.3cm]{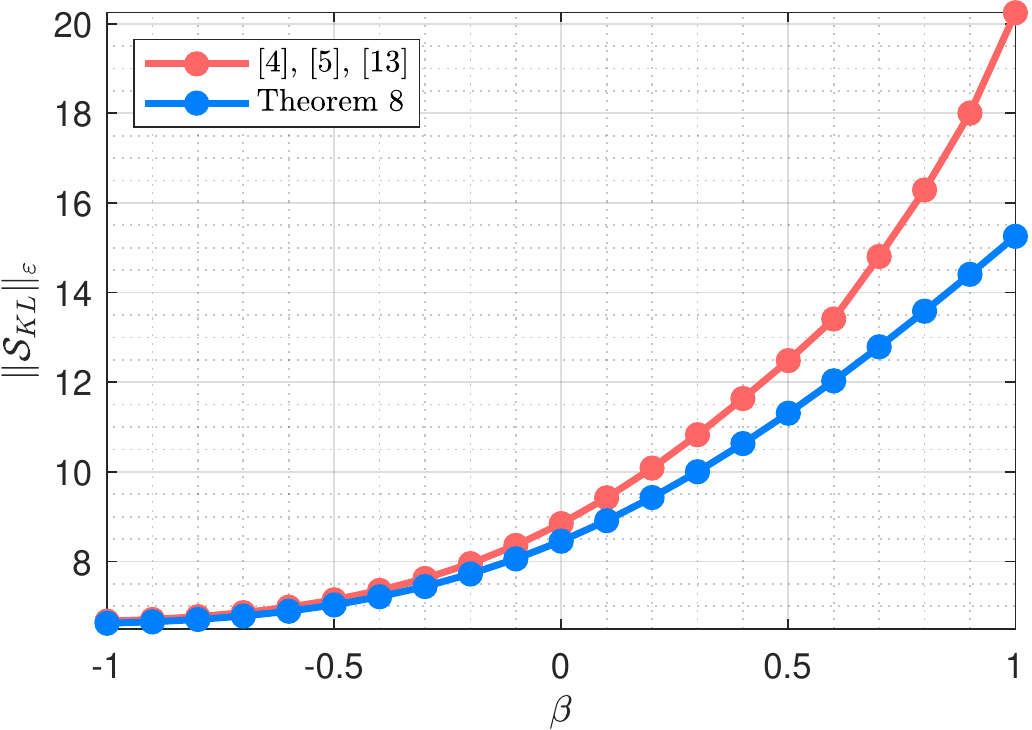}
    \caption{A comparison of the closed-loop system's $\varepsilon$-norm for $\beta \in [-1, 1]$ {{}between the sub-optimal controller \cite{b13-khlebnikov_filter}, \cite{b2-topunov}, \cite{b3-poznyak}} and the optimal controller based on Theorem \ref{thm:KL}.}
    \label{fig3}
\end{figure}

Fig. \ref{fig2} compares $\| \mathcal{S}_{KL} \|_{\varepsilon(\alpha)}$ for {{}$\beta = 0.3$} and $\alpha \in (0, 1)$ between the sub-optimal controller described in {{}\cite{b13-khlebnikov_filter}, \cite{b2-topunov}, \cite{b3-poznyak}} and the optimal controller based on Theorem \ref{thm:KL}. Fig. \ref{fig3} shows a similar comparison of $\| \mathcal{S}_{KL} \|_\varepsilon$ for $\beta \in [-1, 1]$. It is clear that the proposed method provides better results than the previously known one. Table \ref{table2} presents a comparison of the controller parameters obtained with both methods for $\beta = -1$ and $\beta = 1$. 

\begin{table}[h]
    \centering
    \begin{tabular}{cccc}
    \toprule
     \multicolumn{2}{c}{$\beta = -1$} & \multicolumn{2}{c}{$\beta = 1$} \\
     \midrule
         Theorem \ref{thm:KL} & {{} \cite{b13-khlebnikov_filter}, \cite{b2-topunov}, \cite{b3-poznyak}} & Theorem \ref{thm:KL} & {{} \cite{b13-khlebnikov_filter}, \cite{b2-topunov}, \cite{b3-poznyak}} \\
      \cmidrule(r){1-2} \cmidrule(l){3-4} 
      $\hat \alpha = 0.43$ & $\hat \alpha = 0.42$ & $\hat \alpha = 0.82$ & {{}$\hat \alpha \approx 0.4$} \\
      $ \! K \! = \! \begin{bmatrix}
          -0.81 \\ -1.85
      \end{bmatrix}^\top \!$ 
      &
        $\! K \! = \! \begin{bmatrix}
          -0.80 \\ -1.84
      \end{bmatrix}^\top \! $ 
      &
      $\! K \! = \! \begin{bmatrix}
          -3.54 \\ -3.28
      \end{bmatrix}^\top \!$ & 
      {{} $ K \! \approx \! \begin{bmatrix}
          -2.8 \\ -2.6
      \end{bmatrix}^\top \! $ }
      \\ 
       $ \! L  = \! \begin{bmatrix}
          -1.85 \\ -0.81
      \end{bmatrix}^{\phantom \top} \!$ 
      &
       $ \! L  = \! \begin{bmatrix}
          -1.51 \\ -0.55
      \end{bmatrix}^{\phantom \top} \! $ 
      &
      $\! L  = \! \begin{bmatrix}
          -3.28 \\ -3.54
      \end{bmatrix}^{\phantom \top} \!$ & 
      {{} $\! L  \approx \! \begin{bmatrix}
          -1.8 \\ -1.8
      \end{bmatrix}^{\phantom \top} \! $} \\
      \cmidrule(r){1-2} \cmidrule(l){3-4} 
      $\! \|\mathcal{S}_{KL}\|_\varepsilon \!= \! 6.62 \!$ & $\!\|\mathcal{S}_{KL}\|_\varepsilon \! =\!  6.67 \!$ &
      $\!\|\mathcal{S}_{KL}\|_\varepsilon \!= \! 15.3 \!$ & {{} $\!\|\mathcal{S}_{KL}\|_\varepsilon \! \approx \! 20.3 \!$} \\
       \bottomrule
    \end{tabular}
    \caption{Comparison of controller parameters for $\beta = -1, 1$.}
    \label{table2}
\end{table}

{{} Note that for $\beta \ge 0.6$  the sub-optimal controller \cite{b13-khlebnikov_filter}, \cite{b2-topunov}, \cite{b3-poznyak} only provides an approximate solution, as the corresponding LMIs become close to degenerate and the conventional solvers (\texttt{SeDuMi}, \texttt{SDPT3}) encounter numerical problems.}

This example illustrates the contrast between the widely accepted and well-known method and the one proposed in this work, which demonstrates the superiority of the latter for the case of the output-feedback problem.

For the state-feedback and filtering problems the solutions given by Theorems \ref{thm:K} and \ref{thm:L} coincide with the results of \cite{b1-nazin} and \cite{b13-khlebnikov_filter} respectively. However, numerical simulations show that the solution tends to be found much faster with the proposed method, since with modern algorithms for  algebraic Riccati equations tend to be solved more efficiently than the optimization problems with LMI constraints. We also note that, since the separation principle does not hold for the minimization of $\| \mathcal{S}_{KL} \|_\varepsilon$, matrices $K$ and $L$ obtained from \cite{b1-nazin} and \cite{b13-khlebnikov_filter} will not serve as the optimal ones for the output-feedback controller.

\subsection{Comparison with the \texorpdfstring{$\mathcal{H}_2$}{H2}-optimal control}
\label{S5B}
One can observe that the structure of \eqref{ricQ}-\eqref{ricQ-K}, \eqref{ricP}-\eqref{ricP-L} is similar to the solution of the $\mathcal{H}_2$-optimal control problem, as described in \cite{b10-h2}, \cite{b11-mathbook}. However, there are some important differences to highlight.

The proposed approach aims to minimize the value of $\|\mathcal{S}_{KL}\|_{\varepsilon}$, which is particularly relevant for systems subject to bounded disturbances with a finite $\infty$-norm.  In contrast, the $\mathcal{H}_2$-optimal controller minimizes the value of $\|\mathcal{S}_{KL}\|_{\mathcal{H}_2}$, making it more suitable for systems subject to energy-bounded (and therefore essentially decaying) \color{black} disturbances with a finite $2$-norm. Specifically, we have 
    \begin{equation*}
        \|z\|_\infty \le \|\mathcal{S}_{KL}\|_{\varepsilon} \|w\|_\infty, \quad \|z\|_\infty \le \|\mathcal{S}_{KL}\|_{\mathcal{H}_2} \|w\|_2.
    \end{equation*}  

Additionally, the well-known separation principle applies to $\mathcal{H}_2$-control. It states that if $K$ and $L$ are selected as optimal solutions for the state-feedback and filtering problems, respectively, then they will also be the optimal solution for the output-feedback problem. However, this principle does not hold for $\varepsilon$-norm minimization. In general, the matrices $K$ and $L$ that are optimal for the output-feedback problem will differ from those that are optimal for the state-feedback and filtering problems. This is because the minimizing value of $\hat \alpha$ is generally different for each of these three cases.

\section{Discussion}
\label{S6}

In Section \ref{S2}, algebraic equations \eqref{Pinf}, \eqref{Q1} provide ellipsoidal approximations of $\mathcal{R}_\infty$, $\mathcal{O}_1$. We note that these can also be turned into differential matrix equations to provide ellipsoidal approximations of $\mathcal{R}_\infty (T)$, $\mathcal{O}_1 (T)$.

It should be noted that some of the assumptions made in Sections \ref{S3} and \ref{S4} can be relaxed with minimal modifications to the definitions and theorems. For instance, in Section \ref{S4}, we can assume that $(A, B_1)$ is only stabilizable and $(C_2, A)$ is only detectable, which will result in positive \textit{semi}definite solutions $Q_\alpha, P_\alpha \succeq 0$. Moreover, the orthogonality conditions $B_1 D_1^\top = 0$ and $C_2^\top D_2 = 0$ can also be omitted, resulting in more complex {{}versions} of equations \eqref{ricQ}, \eqref{ricP}, but the essence of the result remains unchanged.

\section{Conclusion}
\label{S7}

This paper has investigated the {{}duality relations} between ellipsoidal approximations of {{}easily reachable} and {{}hardly observable} sets for linear systems. By utilizing the duality of the $\varepsilon$-norm, a novel approach for addressing state-feedback and filtering problems has been introduced. The paper's main contribution is the optimal solution for the output-feedback control problem with respect to the $\varepsilon$-norm, outperforming prior results. However, the method still requires one-dimensional iteration to find the optimal solution, and the authors look forward to future research that simplifies this process.

\appendices

\hypertarget{appendix}{}
\section*{Appendix \\ Proofs of Propositions and Theorems} 

\hypersetup{linkcolor={subsectioncolor}}

\subsection{Proofs for Section \ref{S2} -- Analysis}

\hypersetup{linkcolor={blue!70!black}}

\begin{prooftheorem}[\ref{thm:Q1}]
    It is straightforward to check that {{}\eqref{Q1}} is the Lyapunov equation with the solution 
    \begin{equation*}
        Q_\alpha = \int_0^\infty \frac{e^{\alpha t}}{\alpha} e^{A^\top t}C^\top C e^{A t} dt. 
    \end{equation*}
Then for $\mathcal{S}_y$ with $x(0)=x_0$ and $y \in \mathcal{F}^k(T)$ we have
\begin{equation*}
\begin{aligned}
    x_0^\top Q_\alpha x_0 & = \int_0^\infty \frac{e^{\alpha t}}{\alpha} \left| Ce^{At} x_0 \right|^2 dt \\
    & =  \left( \int_0^\infty e^{\alpha t} \left| Ce^{At} x_0 \right|^2 dt \right) \left( \int_0^\infty e^{-\alpha t} dt \right)\\ 
    & \ge \left( \int_0^\infty \left| Ce^{At} x_0 \right| dt \right)^2 \ge \|y\|_1^2,
\end{aligned}
\end{equation*}
where the second to last inequality holds by Cauchy-Schwarz. It follows that $
    x_0^\top Q_\alpha x_0 \le 1 \; \Rightarrow \; \|y\|_1 \le 1. $
\endproof
\end{prooftheorem}


\begin{lemma} 
\label{lemma:1}
\begin{enumerate}
    \item[]
    \item[(i)] If $\tilde Q \succ 0$, $\tilde QA+A^\top \tilde Q \prec 0$, $\tilde Q \succeq C^\top C$, then 
    \begin{equation*}
        \tilde Q \succeq e^{A^\top t}C^\top C e^{At}, \quad \forall t \ge 0.
    \end{equation*}
    \item[(ii)]  If $\tilde P \succ 0$, $A \tilde P+\tilde PA^\top \prec 0$, $\tilde P \succeq B B^\top $, then 
    \begin{equation*}
        \tilde P \succeq e^{At} B B^\top e^{A^\top t}, \quad \forall t \ge 0.
    \end{equation*}
    \vspace{-\baselineskip}
\end{enumerate}
\end{lemma}

\begin{prooflemma}[\ref{lemma:1}]
    From the first two inequalities of (i) we see that $V(x) = x^\top \tilde Q x$ is the Lyapunov function for $\dot x = Ax$. Therefore, $x(t)^\top \tilde Q x(t) \le x(0)^\top \tilde Q x(0)$ for all $t \ge 0$. Hence, for each $x_0$ we have 
    \begin{equation*}
    x_0^\top e^{A^\top t}C^\top C e^{At} x_0 \le x_0^\top e^{A^\top t} \tilde Q e^{At} x_0 \le x_0^\top \tilde Q x_0.    
    \end{equation*}
    If follows that $e^{A^\top t}C^\top C e^{At} \preceq  \tilde Q$ for all $t \ge 0$. To prove (ii) consider the substitution $A \mapsto A^\top$, $C \mapsto B^\top$.
\endproof
\end{prooflemma}

\begin{prooftheorem}[\ref{thm:Qinf}]
    By Lemma \ref{lemma:1} (i) we have
    \begin{equation*}
        |y(t)|^2 = x_0^\top e^{A^\top t}C^\top C e^{At} x_0 \le x_0^\top \tilde Q x_0.
    \end{equation*}
    It follows that $x_0^\top \tilde Q x_0 \le 1 \; \Rightarrow \; \|y\|_\infty \le 1$.
    \endproof
\end{prooftheorem}

\begin{prooftheorem}[\ref{thm:P1}]
      By Lemma \ref{lemma:1} (ii) we have 
      \begin{equation*}
           e^{At} B B^\top e^{A^\top t} \preceq \tilde P,
      \end{equation*}
      which leads to
    \begin{align*}
    &\tilde P^{-1/2}e^{At} B B^\top e^{A^\top t} \tilde P^{-1/2}\preceq  I, \\
    &\lambda_{\text{max}} \left( \tilde P^{-1/2}e^{At} B B^\top e^{A^\top t} \tilde P^{-1/2} \right) \le 1, \\
    &\sigma_{\text{max}} \left( \tilde P^{-1/2}e^{At} B \right) \le 1.
    \end{align*} 
    Therefore, for $u \in \mathcal{F}^m(T)$ and $t \le T$ we have
    \begin{align*}
         |\tilde P^{-1/2}x(t)| &= \left| \int_0^t \tilde P^{-1/2}e^{A(t-\tau)}Bu(\tau) d\tau \right|  \\
         & \le  \int_0^t \left| \tilde P^{-1/2}e^{A(t-\tau)}Bu(\tau) \right| d\tau \\
         & \le \int_0^t \sigma_{\text{max}} \left( \tilde P^{-1/2}e^{A(t-\tau)}B \right) \left| u(\tau) \right|  d\tau \\
         & \le \sup_{t \ge 0} \left( \sigma_{\text{max}} \left( \tilde P^{-1/2}e^{At}B \right) \right) \int_0^t  \left| u(\tau) \right|  d\tau \\
         &  \le \int_0^t  \left| u(\tau) \right|  d\tau \le  \|u\|_1.
    \end{align*} 
    Hence, \\
    \begin{equation*}
         \|u\|_1 \le 1 \; \Rightarrow \; |\tilde P^{-1/2}x(t)|^2 = x(t)^\top \tilde P^{-1} x(t) \le 1.
    \end{equation*}
This completes the proof.
\endproof
\end{prooftheorem}

\begin{proofproposition}[\ref{pr:trace}]
    For $\alpha \in (0, -2r)$ the system \eqref{true_system} is stable. The claim follows from Remark \ref{remark1} and the fact about controllability and observability Gramians (see, e.g., \cite{b17-h2_norm}).
    \endproof
\end{proofproposition}

\begin{proofproposition}[\ref{pr:max}]
Using Schur complement properties, one can check that if $\tilde P \succ 0$ is a solution to 
\begin{equation}
\label{P1full}
     A \tilde P + \tilde P A^\top \prec 0, \quad \tilde P \succeq BB^\top, \quad C\tilde P C^\top \preceq \lambda I \phantom{.}
\end{equation}
for some $\lambda > 0$, then $\tilde Q = \lambda \tilde P^{-1} \succ 0 $ is a solution to
\begin{equation}
\label{Q1full}
     \tilde Q  A + A^\top  \tilde Q \prec 0, \quad  \tilde Q \succeq C^\top C, \quad B^\top \tilde Q B \preceq \lambda I.
\end{equation}
Similarly, if $\tilde Q \succ 0$ is a solution to \eqref{Q1full}, then $\tilde P = \lambda \tilde Q^{-1} \succ 0$ is a solution to \eqref{P1full}. So, for each $\lambda$, for which there is a solution $\tilde P$, there is also a solution $\tilde Q$, and vice versa. Therefore, the minimum values of $\lambda$ for which  \eqref{P1full} and \eqref{Q1full} are feasible coincide. By the last inequalities of \eqref{P1full} and \eqref{Q1full}, they are exactly the maximum eigenvalues of $C \tilde P C^\top$ and $B^\top \tilde Q B$, hence the claim. 
\endproof
\end{proofproposition}

\begin{proofproposition}[\ref{pr:Pinequality}]
    It is known that if $x_\ast^\top \mathcal{P}^{-1} x_\ast^{\phantom{\top}} \le 1$ and $y_\ast = Cx_\ast$, then $y_\ast^\top (C \mathcal{P} C^\top)^{-1} y_\ast^{\phantom{\top}} \le 1$. Therefore,
    \begin{equation*}
    \begin{aligned}
    \| \mathcal{S} \|_{\infty, p} & = \sup_{T \ge 0} \: \max_{\|u\|_p \le 1} \, \|y\|_\infty 
     = \max_{x_\ast \in \mathcal{R}_p} |Cx_\ast| \\
    & \le \max_{x_\ast^\top \mathcal{P}^{-1} x_\ast^{\phantom{\top}} \le 1} |Cx_\ast| 
     \: = \max_{y_\ast^\top (C \mathcal{P} C^\top)^{-1} y_\ast^{\phantom{\top}} \le 1} |y_\ast| \\ &  
    = \sqrt{\lambda_{\text{max}} (C \mathcal{P} C^\top)},
    \end{aligned}      
    \end{equation*}
hence the first inequality of the claim. The second inequality holds because of the general properties of $\operatorname{trace}$. 
\endproof
\end{proofproposition}

\begin{proofproposition}[\ref{pr:Qinequality}]
    Note that $\mathcal{S}$ with $x(0)=0$ and $u(t) = u_0 \delta (t)$ is equivalent to $\mathcal{S}$ with $x(0)=Bu_0$ and $u=0$, i.e. an impulse provides the system with initial conditions. If we consider only $x_0$ of the form $x_0 = B u_0$, then
    \begin{equation*}
    \begin{aligned}
    \| \mathcal{S} \|^{-1}_{q, i} & = \inf_{T \ge 0} \: \min_{\|y\|_q \ge 1} \, |u_0| = \inf_{x_0 \not\in \mathcal{O}_q} |u_0| \\
    & \ge \min_{x_0^\top \mathcal{Q} x_0^{\phantom{\top}} \ge 1} |u_0| \quad   = \min_{u_0^\top (B^\top \mathcal{Q} B) u_0^{\phantom{\top}} \ge 1} |u_0| \\
    & = \frac{1}{\sqrt{\lambda_{\text{max}} (B^\top \mathcal{Q} B)}},
    \end{aligned}      
    \end{equation*}
hence the first inequality of the claim. The second inequality holds because of the general properties of $\operatorname{trace}$.
\endproof
\end{proofproposition}

\begin{prooftheorem}[\ref{thm:norms}]
    Follows from Propositions \ref{pr:Pinequality} and \ref{pr:Qinequality}.
    \endproof
\end{prooftheorem}

\hypersetup{linkcolor={subsectioncolor}}

\subsection{Proofs for Section \ref{S3} -- State-Feedback and Filtering}

\hypersetup{linkcolor={blue!70!black}}

\begin{proofproposition}[\ref{pr:Q-solvable}]
    Note that \eqref{ricQ} is a standard Riccati equation for the system \eqref{true_system}, so it has the unique positive definite solution $Q_\alpha \succ 0$, whenever
    \begin{align*}
        & \Big( A+\frac{\alpha}{2}I, \, \frac{1}{\sqrt{\alpha}}B \Big) \text{ is stabilizable}, \\
        &\Big(  \frac{1}{\sqrt{\alpha}}C, \, A+\frac{\alpha}{2}I \Big) \text{ is observable}. 
    \end{align*}
    Using the Popov-Belevitch-Hautus test, it is straightforward to show that both of these hold for all $\alpha>0$, if $(A,B)$ is stabilizable and $(C,A)$ is observable. Hence, for each $\alpha > 0$ there exists the corresponding $Q_\alpha \succ 0$.
    \endproof
\end{proofproposition}

\begin{prooftheorem}[\ref{thm:K}]
   Consider the closed-loop system 
    \begin{equation*}
       \mathcal{S}_K : \;  \left\{ 
        \begin{aligned}
            &\dot x = (A+BK)x + B_w w, \\
            &z = (C+DK)x,
        \end{aligned}
        \right.
    \end{equation*}
    and apply equation \eqref{Q1} to obtain 
    \begin{equation*} \begin{multlined}[c][.9\displaywidth]
        Q_\alpha(A+BK)+(A+BK)^\top Q_\alpha  \\ + \alpha Q_\alpha + \frac{1}{\alpha} (C+DK)^\top  (C+DK) =  0. 
    \end{multlined} \end{equation*} 
    With the assumption $C^\top D = 0$, we get
    \begin{equation*} \begin{multlined}[c][.9\displaywidth]
        Q_\alpha A + A^\top Q_\alpha  +Q_\alpha BK + K^\top B^\top Q_\alpha \\ + \alpha Q_\alpha + \frac{1}{\alpha}C^\top C  + \frac{1}{\alpha}K^\top D^\top D K  =  0. 
    \end{multlined} \end{equation*} 
    Given that $D^\top D$ is invertible, completing the square gives
    \begin{align*}
    &\frac{1}{\alpha} \! \left(K+\alpha (D^\top \! D)^{-1}B^\top  Q_\alpha\right)^\top \! D^\top \! D \left(K+\alpha (D^\top \! D)^{-1}B^\top  Q_\alpha\right) \\
    & \quad \quad \quad \quad \quad \quad \quad \; \; \; + Q_\alpha A + A^\top Q_\alpha + \alpha Q_\alpha \\ 
    & \quad \quad \quad \quad \quad \quad \quad \; \; \;  - \alpha Q_\alpha B(D^\top \! D)^{-1} B^\top Q_\alpha + \frac{1}{\alpha} C^\top C  =0. 
    \end{align*}
    The first term is a square, hence positive semidefinite, so the rest is negative semidefinite, i.e.
    \begin{equation*}
      \begin{multlined}[c][.85\displaywidth]
    Q_\alpha A + A^\top Q_\alpha + \alpha Q_\alpha \\ 
    - \alpha Q_\alpha B(D^\top D)^{-1} B^\top Q_\alpha + \frac{1}{\alpha} C^\top C  \preceq 0.
      \end{multlined}  
    \end{equation*}
    But note that if $X, Y \succ 0$ are such that
    \begin{alignat*}{6}
        & X A + && A^\top X + && \alpha X - && \alpha X B(D^\top D)^{-1} B^\top X + &&\alpha^{-1} C^\top C  =  0, \\
        &Y A + &&A^\top Y + &&\alpha Y - &&\alpha Y B(D^\top D)^{-1} B^\top Y + &&\alpha^{-1} C^\top C  \preceq 0,
    \end{alignat*}
    then $X \preceq Y$ (see \cite{b15-riccati_comparison}), and
    \begin{equation*}
        \operatorname{trace} (B_w^\top X B_w^{\phantom \top}) \le \operatorname{trace} (B_w^\top Y B_w^{\phantom \top}).
    \end{equation*}
    It means that the minimum value of the $\varepsilon(\alpha)$-norm is achieved when \eqref{ricQ} holds, so \eqref{ricQ-K} holds. Equation \eqref{ricQ} admits the unique positive definite solution $Q_\alpha$ by Proposition \ref{pr:Q-solvable}. From the general theory of Riccati equations (see, e.g., \cite{b15-riccati_comparison}, \cite{b14-riccati_equation}) the corresponding $K$ makes the matrix $A+\frac{\alpha}{2} I+BK$ stable, hence $A+BK$ is stable. Then by the definition of the $\varepsilon (\alpha)$-norm we have
    \begin{equation*}
         \| \mathcal{S}_K \|_{\varepsilon(\alpha)}^2 =  \operatorname{trace} (B_w^\top Q_\alpha B_w^{\phantom \top}).
    \end{equation*}
    This is exactly the claim of Theorem \ref{thm:K}. 
    \endproof
\end{prooftheorem}

\begin{proofproposition}[\ref{pr:notzero}]
    For small $\alpha$ we have $\alpha Q_\alpha \approx \hat Q$, where $\hat Q$ is defined by
    \begin{equation*}
        \hat Q A + A^\top \hat Q - \hat Q B(D^\top D)^{-1} B^\top \hat Q +C^\top C = 0.
    \end{equation*}
    Therefore, if $\alpha \approx 0$, then 
    \begin{equation*}
        \operatorname{trace} (B_w^{\top} Q_{\alpha} B_w^{\phantom \top}) \approx \frac{1}{\alpha} \operatorname{trace} (B_w^{\top} \hat Q B_w^{\phantom \top}), 
    \end{equation*}
    which is a decaying function of $\alpha$. Hence, the infimum is achieved away from zero.
    \endproof
\end{proofproposition}

\begin{proofproposition}[\ref{pr:P-solvable}]  
    Dual to Proposition \ref{pr:Q-solvable}.
    \endproof
\end{proofproposition}

\begin{prooftheorem}[\ref{thm:L}]
   Consider the closed-loop system 
    \begin{equation*}
       \mathcal{S}_L : \;  \left\{ 
        \begin{aligned}
            &\dot e = (A+LC)e + (B+LD)w, \\
            &z = C_z e,
        \end{aligned}
        \right.
    \end{equation*}
    where $e=x-\hat x$, and apply equation \eqref{Pinf} to obtain 
    \begin{equation*} \begin{multlined}[c][.9\displaywidth]
        (A+LC) P_\alpha + P_\alpha (A+LC)^\top    \\ + \alpha P_\alpha + \frac{1}{\alpha} (B+LD)(B+LD)^\top =  0. 
    \end{multlined} \end{equation*} 
    The rest of the proof is dual to the one of Theorem \ref{thm:K}. 
    \endproof
\end{prooftheorem}

\begin{proofproposition}[\ref{pr:notzero2}]
    Dual to Proposition \ref{pr:notzero}.
    \endproof
\end{proofproposition}

\hypersetup{linkcolor={subsectioncolor}}

\subsection{Proof of Theorem \ref{thm:KL} of Section \ref{S4} -- Main result on optimal output-feedback control with respect to \texorpdfstring{$\varepsilon$}{ε}-norm}

\hypersetup{linkcolor={blue!70!black}}

Let $\mathcal{S}_1$ be a system with input $w$ and output $z_1$, let $\mathcal{S}_2$ be a system with input $w$ and output $z_2$. Then we define the \textit{sum} $\mathcal{S}_1 + \mathcal{S}_2$ as the system with input $w$ and output $z = z_1+z_2$.


\begin{figure}[H]
    \centering
\begin{tikzpicture}
\node [] (start0) at (0,0) {$w$};
\node [draw,
    minimum width=1.6cm,
    minimum height=0.8cm,
    right = 0.4 cm of start0,
]  (s) {$\mathcal{S}_{1}+\mathcal{S}_2$};
\node [right = 0.4 cm of s] (finish0)  {$z$};
\draw[-latex] (start0) -- (s.west);
\draw[-latex] (s.east) -- (finish0);
\node [right = 0.15cm of finish0] (equal)  {\LARGE $\vcentcolon =$};
\node [right = 0.15cm of equal] (before_start) {$w$};
\node [right = 0.15cm of before_start] (start) {};
\node [draw,
    minimum width=1cm,
    minimum height=0.8cm,
    above right = 0.08 cm and 0.4 cm of start, 
]  (s1) {$\mathcal{S}_1$};
\node [draw,
    minimum width=1cm,
    minimum height=0.8cm,
    below right = 0.08 cm and 0.4cm of start,
]  (s2) {$\mathcal{S}_2$};
\node[draw,
    circle,
    minimum size=0.4cm,
    label = center:$+$,
    above right = 0.083cm and 0.4cm of s2,
] (sum) {};
\node [right = 0.25 cm of sum] (finish) {$z$};
\draw[-latex] (start.center) |- (s1.west);
\draw[-latex] (start.center) |- (s2.west);
\draw[-latex] (s1.east) -| (sum.north)
        node[pos = 0.21,above]{$z_1$};
\draw[-latex] (s2.east) -| (sum.south)
        node[pos = 0.21, above]{$z_2$};
\draw[-latex] (sum.east) -- (finish);
\draw (before_start.east) -- (start.center);
\end{tikzpicture}
\end{figure}

    Let $\mathcal{S}_1$ be a system with input $w$ and output $v$, let $\mathcal{S}_2$ be a system with input $v$ and output $z$. Then we define the \textit{product} $\mathcal{S}_2 \mathcal{S}_1$ as the system with input $w$ and output $z$.

\begin{figure}[H]
    \centering
\begin{tikzpicture}
\node [right = 0.4cm of equal] (start)  {$z$};
\node [draw,
    minimum width=1.2cm,
    minimum height=0.8cm,
    right = 0.4 cm of start,
]  (s) {$\mathcal{S}_2 \mathcal{S}_1$};
\node [right = 0.4 cm of s] (finish)  {$w$};
\node [right = 0.15cm of finish] (equal)  {\LARGE $\vcentcolon =$};
\node [right = 0.15cm of equal] (start0)  {$z$};
\node [draw,
    minimum width=1cm,
    minimum height=0.8cm,
    right = 0.4 cm of start0,
]  (s2) {$\mathcal{S}_2$};
\node [draw,
    minimum width=1cm,
    minimum height=0.8cm,
    right = 0.6 cm of s2,
]  (s1) {$\mathcal{S}_1$};
\node [right = 0.4 cm of s1] (finish0)  {$w$};
\draw[-latex] (s.west) -- (start);
\draw[-latex] (finish) -- (s.east);
\draw[-latex] (s2.west) -- (start0);
\draw[-latex] (s1.west) -- (s2.east)
        node[pos = 0.45,above]{$v$};
\draw[-latex] (finish0) -- (s1.east);
\end{tikzpicture}
\end{figure}

We {{}are ready} to state the following lemmas. 

\begin{lemma}
\label{lemma:sum-norm}
If $\mathcal{S}_1+\mathcal{S}_2$ is well-defined, then
\begin{equation*}
    \| \mathcal{S}_1 + \mathcal{S}_2\|_{\varepsilon(\alpha)}^2 = 
    \| \mathcal{S}_1 \|_{\varepsilon(\alpha)}^2 + \| \mathcal{S}_2\|_{\varepsilon(\alpha)}^2.
\end{equation*}
\vspace{-\baselineskip}
\end{lemma}

\begin{prooflemma}[\ref{lemma:sum-norm}] Let $\mathcal{S}_1$ and $\mathcal{S}_2$ be given as
    \begin{equation*}
    \mathcal{S}_i : \; \left\{ 
    \begin{aligned}
        &\dot x_i = A_ix_i + B_iw, \\
        &z_i = C_i x_i,
    \end{aligned}
    \right.
    \quad i = 1,2,
\end{equation*}
and let $P_i = P_{\alpha, i}$ and $Q_i = Q_{\alpha, i}$ be the solutions of the equations of types \eqref{Pinf} and \eqref{Q1} for some $\alpha > 0$, so that 
\begin{equation*}
     \| \mathcal{S}_i\|_{\varepsilon(\alpha)}^2 =  \operatorname{trace} (C_i P_i C^\top_i) =  \operatorname{trace} (B^\top_i Q_i B_i), \quad i = 1, 2.
\end{equation*}
Then it is straightforward to check, that
\begin{equation*}
    P_\alpha  = \begin{bmatrix}
        P_{1} & 0 \\ 0 & P_2
    \end{bmatrix}, \quad 
    Q_\alpha = \begin{bmatrix}
        Q_1 & 0 \\ 0 & Q_2
    \end{bmatrix}
\end{equation*}
are the solutions of \eqref{Pinf} and \eqref{Q1}  for $\mathcal{S} = \mathcal{S}_1 + \mathcal{S}_2$ with
   \begin{equation*}
      A = \begin{bmatrix}
        A_1 & 0 \\ 
        0 & A_2
    \end{bmatrix}, \; 
    B = \begin{bmatrix}
        B_1 \\ B_2
    \end{bmatrix}, \;
    C = \begin{bmatrix}
        C_1 & C_2
    \end{bmatrix}, 
\end{equation*}
and $\| \mathcal{S}\|_{\varepsilon(\alpha)}^2 = \sum \operatorname{trace} (C_i P_i C^\top_i)  = \sum \operatorname{trace} (B^\top_i Q_i B_i)$.
\vspace{0.2\baselineskip}
\endproof
\end{prooflemma}

\begin{lemma} 
\label{lemma:product-norm}
    (i) Let $\mathcal{S}_1$, $\mathcal{S}_1^\prime$ and  $\mathcal{S}_2^*$ be given as
    \begin{equation*}
    \mathcal{S}_1 : \; \left\{ 
    \begin{aligned}
        &\dot x_1 = A_1x_1 + B_1w, \\
        &z_1 = C_1 x_1,
    \end{aligned}
    \right.
    \quad 
    \mathcal{S}_1^\prime : \; \left\{ 
    \begin{aligned}
        &\dot x^\prime= A_1 x^\prime + B_1w, \\
        &z ^\prime = D_2 C_1 x^\prime,
    \end{aligned}
    \right.
    \end{equation*}
    \begin{equation*}
    \mathcal{S}_2^* : \; \left\{ 
    \begin{aligned}
        &\dot x_2 = (A_2 + B_2 K) x_2 + B_2 z_1, \\
        &z_2 = (C_2 + D_2 K) x_2 + D_2 z_1,
    \end{aligned}
    \right.
    \end{equation*}
where $K$ is calculated from \eqref{ricQ}, \eqref{ricQ-K} with $(A,B,C,D) = (A_2,B_2,C_2,D_2)$, and $C_2^\top D_2 =0$. Then 
\begin{equation*}
    \|\mathcal{S}_2^* \mathcal{S}_1^{\phantom \prime}\|_{\varepsilon(\alpha)} = \|\mathcal{S}_1^\prime \|_{\varepsilon(\alpha)}.
\end{equation*}
(ii) Let $\mathcal{S}_2$, $\mathcal{S}_2^\prime$ and  $\mathcal{S}_1^*$ be given as
    \begin{equation*}
    \mathcal{S}_2 : \; \left\{ 
    \begin{aligned}
        &\dot x_2 = A_2x_2 + B_2 z_1, \\
        &z_2 = C_2 x_2,
    \end{aligned}
    \right.
    \quad 
    \mathcal{S}_2^\prime : \; \left\{ 
    \begin{aligned}
        &\dot x^\prime= A_2 x^\prime + B_2 D_1 w, \\
        &z ^\prime = C_2 x^\prime,
    \end{aligned}
    \right.
    \end{equation*}
    \begin{equation*}
    \mathcal{S}_1^* : \; \left\{ 
    \begin{aligned}
        &\dot x_1 = (A_1 + LC_1) x_1 + (B_1+L D_1)w, \\
        &z_1 = C_1 x_1 + D_1 w,
    \end{aligned}
    \right.
    \end{equation*}
where $L$ is calculated from \eqref{ricP}, \eqref{ricP-L} with $(A,B,C,D) = (A_1,B_1,C_1,D_1)$, and $B_1 D_1^\top =0$. Then 
\begin{equation*}
    \|\mathcal{S}_2^{\phantom \prime} \mathcal{S}_1^*\|_{\varepsilon(\alpha)} = \|\mathcal{S}_2^\prime \|_{\varepsilon(\alpha)}.
\end{equation*}
\vspace{-\baselineskip}
\end{lemma}

\begin{prooflemma}[\ref{lemma:product-norm}]
    (i) Let $Q_\alpha \succ 0$ and $K$ be the solutions of \eqref{ricQ}, \eqref{ricQ-K} with $(A,B,C,D)=(A_2,B_2,C_2,D_2)$ for some $\alpha >0$.
    Let $Q_\alpha^\prime \succ 0$ be the solution of 
    \begin{equation*}
        Q_\alpha^\prime A_1 + A_1^\top Q_\alpha^\prime + \alpha Q_\alpha ^\prime \\ + \frac{1}{\alpha}C_1^\top D_2^\top D_2 C_1 = 0,
    \end{equation*}
    which gives $\|\mathcal{S}_1^\prime \|_{\varepsilon(\alpha)}^2 = \operatorname{trace} (B_1^\top Q_\alpha^\prime B_1)$.    
    Then
    \begin{equation*}
        \bar Q_\alpha = \begin{bmatrix}
            Q_\alpha^\prime & 0 \\ 0 & Q_\alpha
        \end{bmatrix}
    \end{equation*}
    is the solution to \eqref{Q1} for the system $\bar {\mathcal{S}} = \mathcal{S}_2^* \mathcal{S}_1^{\phantom \prime}$ with 
    \begin{align*}
        \bar A & = \begin{bmatrix}
        A_1 &  0\\ 
        B_2 C_1 & A_2 + B_2 K
    \end{bmatrix}, \quad 
    \bar B = \begin{bmatrix}
        B_1 \\ 0
    \end{bmatrix}, \\
   \bar C & = \begin{bmatrix}
         D_2 C_1 & C_2 + D_2 K
    \end{bmatrix}. 
    \end{align*}
Therefore,
\begin{equation*}
  \|\mathcal{S}_2^* \mathcal{S}_1^{\phantom \prime}\|_{\varepsilon(\alpha)}^2 = \operatorname{trace} (\bar B^\top \bar Q_\alpha \bar B)= \operatorname{trace} (B_1^\top Q_\alpha^\prime B_1).   
\end{equation*}

(ii) Let $P_\alpha \succ 0$ and $L$ be the solutions of \eqref{ricP}, \eqref{ricP-L} with $(A,B,C,D)=(A_1,B_1,C_1,D_1)$ for some $\alpha >0$.    Let $P_\alpha^\prime \succ 0$ be the solution of 
    \begin{equation*}
        A_2 P_\alpha^\prime  + P_\alpha^\prime  A_2^\top + \alpha P_\alpha^\prime \\ + \frac{1}{\alpha}B_2 D_1 D_1^\top B_2^\top = 0,
    \end{equation*}
    which gives $\|\mathcal{S}_2^\prime \|_{\varepsilon(\alpha)}^2 = \operatorname{trace} (C_2 P_\alpha^\prime C_2^\top)$.    
    Then
    \begin{equation*}
        \bar P_\alpha = \begin{bmatrix}
            P_\alpha & 0 \\ 0 & P_\alpha^\prime
        \end{bmatrix}
    \end{equation*}
    is the solution to \eqref{Pinf} for the system $\bar {\mathcal{S}} = \mathcal{S}_2^{\phantom \prime} \mathcal{S}_1^*$ with 
    \begin{equation*}  
        \bar A = \begin{bmatrix}
        A_1+LC_1 &  0\\ 
        B_2 C_1 & A_2
    \end{bmatrix}, \;
    \bar B = \begin{bmatrix}
        B_1 + LD_1 \\ B_2 D_1
    \end{bmatrix}, \;
   \bar C = \begin{bmatrix}
         0 & C_2
    \end{bmatrix}. 
 \end{equation*}
Therefore,
\begin{equation*}
  \|\mathcal{S}_2^{\phantom \prime} \mathcal{S}_1^*\|_{\varepsilon(\alpha)}^2 = \operatorname{trace} (\bar C \bar P_\alpha \bar C^\top)= \operatorname{trace} (C_2  P_\alpha^\prime C_2^\top).   
\end{equation*}
This completes the proof of Lemma \ref{lemma:product-norm}.
\endproof
\end{prooflemma}

\begin{prooftheorem} [\ref{thm:KL} (Main result)]
    Consider two following equivalent state-space representations of  $\mathcal{S}_{KL}$:
    \begin{alignat*}{3}
        \begin{bmatrix}
            \dot{x} \\ \dot{e}
        \end{bmatrix} 
        & = 
        &&
        \begin{bmatrix}
            A+B_2 K & -B_2 K \\ 0 & A+LC_1
        \end{bmatrix}
        &&\begin{bmatrix}
            x \\ e
        \end{bmatrix} +
        \begin{bmatrix}
            B_1 \\ B_1 + LD_1
        \end{bmatrix}w , \\
        z 
        & = 
        &&\begin{bmatrix}
            \; C_2 + D_2 K &  - D_2 K \phantom{;}
        \end{bmatrix}
        && \begin{bmatrix}
            x \\ e
        \end{bmatrix},
        \end{alignat*}
and 
    \begin{alignat*}{3}
        \begin{bmatrix}
            \dot{\hat x} \\ \dot{e}
        \end{bmatrix} 
        & = 
        &&
        \begin{bmatrix}
            A+B_2 K & -LC_1 \\ 0 & A+LC_1
        \end{bmatrix}
        &&\begin{bmatrix}
            \hat x \\ e
        \end{bmatrix} +
        \begin{bmatrix}
            -LD_1 \\ B_1 + LD_1
        \end{bmatrix}w , \\
        z 
        & = 
        &&\begin{bmatrix}
            \; C_2 + D_2 K &  C_2 \phantom{;}
        \end{bmatrix}
        && \begin{bmatrix}
            \hat x \\ e
        \end{bmatrix},
        \end{alignat*}
where $e = x - \hat x$. 
Define the following systems:
\begin{alignat*}{2}
       &\mathcal{S}_K : \;  && \left\{ 
        \begin{aligned}
           & \dot x_1  = (A+B_2 K)x_1 + B_1 w, \\
           & v_1 = (C_2+D_2K)x_1,
        \end{aligned}
        \right. \\      
    &\mathcal{S}_L : \;  && \left\{ 
        \begin{aligned}
           & \dot e  = (A+LC_1)e + (B_1+LD_1) w, \\
           & z_1 = C_2 e,
        \end{aligned}
        \right. \\    
    &\mathcal{T}_K^* : \;  && \left\{ 
    \begin{aligned}
       & \dot x_2  = (A+B_2 K)x_2 + B_2 g, \\
       & v_2 = (C_2+D_2K)x_2 + D_2 g,
    \end{aligned}
    \right. \\
    &\mathcal{T}_L^* : \;  && \left\{ 
        \begin{aligned}
           & \dot e  = (A+LC_1)e + (B_1+LD_1) w, \\
           & h = C_1 e + D_1 w,
        \end{aligned}
        \right. \\
    &\mathcal{T}_{LK}^{\phantom *} : \;  && \left\{ 
            \begin{aligned}
               & \dot e  = (A+LC_1)e + (B_1+LD_1) w, \\
               & g = -Ke,
            \end{aligned}
            \right. \\
        &\mathcal{T}_{KL}^{\phantom *} : \;  && \left\{ 
    \begin{aligned}
       & \dot {\hat x}  = (A+B_2 K) \hat x -Lh, \\
       & z_2 = (C_2+D_2K)\hat x.
    \end{aligned}
    \right. 
\end{alignat*}

Note that the definitions of $\mathcal{S}_K$ and $\mathcal{S}_L$ coincide with the ones from Section \ref{S3} (see the proofs of Theorems \ref{thm:K} and \ref{thm:L}). Observe that $z=z_1+z_2=v_1+v_2$ and
\begin{equation*}
    \mathcal{S}_{KL} = \mathcal{S}_K + \mathcal{T}_K^* \mathcal{T}_{LK}^{\phantom *} = \mathcal{S}_L + \mathcal{T}_{KL}^{\phantom *} \mathcal{T}_{L}^*.
\end{equation*} %
\vspace{-\baselineskip}
\begin{figure}[H]
    \centering
\begin{tikzpicture}
\node [] (start) at (0,0) {};
\node [left = 0.15cm of start] (before_start) {$w$};
\node [draw,
    minimum width=1cm,
    minimum height=0.8cm,
    above right = 0.08 cm and 1.2 cm of start, 
]  (s1) {$\mathcal{S}_K
$};
\node [draw,
    minimum width=1cm,
    minimum height=0.8cm,
    below right = 0.08 cm and 0.4cm of start,
]  (s2) {$\mathcal{T}_{LK}^{\phantom *}$};
\node [draw,
    minimum width=1cm,
    minimum height=0.8cm,
    right = 0.6cm of s2,
]  (s3) {$\mathcal{T}_K^\ast$};
\node[draw,
    circle,
    minimum size=0.4cm,
    label = center:$+$,
    above right = 0.083cm and 0.4cm of s3,
] (sum) {};
\node [right = 0.25 cm of sum] (finish) {$z$};
\draw[-latex] (start.center) |- (s1.west);
\draw[-latex] (start.center) |- (s2.west);
\draw[-latex] (s1.east) -| (sum.north)
        node[pos = 0.095,above]{$v_1$};
\draw[-latex] (s2.east) -- (s3.west)
        node[pos = 0.45,above]{$g$};
\draw[-latex] (s3.east) -| (sum.south)
        node[pos = 0.21, above]{$v_2$};
\draw[-latex] (sum.east) -- (finish);
\draw (before_start.east) -- (start.center);

\node [] (start2) at (0,-2.6) {};
\node [left = 0.15cm of start2] (before_start2) {$w$};
\node [draw,
    minimum width=1cm,
    minimum height=0.8cm,
    above right = 0.08 cm and 1.2 cm of start2, 
]  (s4) {$\mathcal{S}_L
$};
\node [draw,
    minimum width=1cm,
    minimum height=0.8cm,
    below right = 0.08 cm and 0.4cm of start2,
]  (s5) {$\mathcal{T}_{L}^\ast$};
\node [draw,
    minimum width=1cm,
    minimum height=0.8cm,
    right = 0.6cm of s5,
]  (s6) {$\mathcal{T}_{KL}^{\phantom *}$};
\node[draw,
    circle,
    minimum size=0.4cm,
    label = center:$+$,
    above right = 0.083cm and 0.4cm of s6,
] (sum2) {};
\node [right = 0.25 cm of sum2] (finish2) {$z$};
\draw[-latex] (start2.center) |- (s4.west);
\draw[-latex] (start2.center) |- (s5.west);
\draw[-latex] (s4.east) -| (sum2.north)
        node[pos = 0.095,above]{$z_1$};
\draw[-latex] (s5.east) -- (s6.west)
        node[pos = 0.45,above]{$h$};
\draw[-latex] (s6.east) -| (sum2.south)
        node[pos = 0.21, above]{$z_2$};
\draw[-latex] (sum2.east) -- (finish2);
\draw (before_start2.east) -- (start2.center);

\node [] (start0) at (-4,-1.3) {$w$};
\node [draw,
    minimum width=1cm,
    minimum height=0.8cm,
    right = 0.4 cm of start0,
]  (s) {$\mathcal{S}_{KL}$};
\node [right = 0.4 cm of s] (finish0)  {$z$};
\draw[-latex] (start0) -- (s.west);
\draw[-latex] (s.east) -- (finish0);
\node [rotate = 30] (u1) at (-1.3,-0.4) {\LARGE $=$};
\node [rotate = -30] (l1) at (-1.3,-2.2) {\LARGE $=$};
\end{tikzpicture}
\end{figure}

\noindent Use Lemmas \ref{lemma:sum-norm} and \ref{lemma:product-norm} to find that 
\begin{alignat*}{3}
    \| \mathcal{S}_{KL} \|_{\varepsilon(\alpha)}^2 
    & = \| \mathcal{S}_K \|_{\varepsilon(\alpha)}^2 && + \| \mathcal{T}_K^* \mathcal{T}_{LK}^{\phantom *} \|_{\varepsilon(\alpha)}^2 \\
    & = \| \mathcal{S}_K \|_{\varepsilon(\alpha)}^2 &&+ \| \mathcal{S}_L^\prime \|_{\varepsilon(\alpha)}^2, \\
    \| \mathcal{S}_{KL} \|_{\varepsilon(\alpha)}^2 
    & = \| \mathcal{S}_L \|_{\varepsilon(\alpha)}^2 &&+ \| \mathcal{T}_{KL}^{\phantom *} \mathcal{T}_{L}^* \|_{\varepsilon(\alpha)}^2 \\ 
    & = \| \mathcal{S}_L \|_{\varepsilon(\alpha)}^2 &&+ \| \mathcal{S}_K^\prime \|_{\varepsilon(\alpha)}^2,
\end{alignat*}
where 
    \begin{alignat*}{2}
    &\mathcal{S}_{L}^\prime : \;  && \left\{ 
            \begin{aligned}
               & \dot e  = (A+LC_1)e + (B_1+LD_1) w, \\
               & v^\prime = -D_2 Ke,
            \end{aligned}
            \right. \\
        &\mathcal{S}_{K}^\prime : \;  && \left\{ 
    \begin{aligned}
       & \dot {x}^\prime  = (A+B_2 K) x^\prime -L D_1 w, \\
       & z^\prime = (C_2+D_2K)x^\prime.
    \end{aligned}
    \right. 
\end{alignat*}

By Theorems \ref{thm:K} and \ref{thm:L}, if we take
\begin{equation*}
    K = - \alpha (D_2^\top D_2)^{-1}B_2^\top Q_\alpha, \quad
    L = - \alpha P_\alpha C_1^\top (D_1 D_1^\top)^{-1}, 
\end{equation*}
then both $\| \mathcal{S}_K \|_{\varepsilon(\alpha)}^2$ and $\| \mathcal{S}_L \|_{\varepsilon(\alpha)}^2$ attain minimum values:
\begin{alignat*}{2}
    & \min_K \| \mathcal{S}_K \|_{\varepsilon(\alpha)}^2 &&= \operatorname{trace} (B_1^\top Q_\alpha B_1), \\
    & \min_L \| \mathcal{S}_L \|_{\varepsilon(\alpha)}^2 &&= \operatorname{trace} (C_2 P_\alpha C_2^\top). 
 \end{alignat*}

Since $\mathcal{S}_K^\prime$ and $\mathcal{S}_L^\prime$ differ from $\mathcal{S}_K$ and $\mathcal{S}_L$ only in input and output matrices, they correspond to the same Riccati equations.
 Therefore, the minimum values for $\| \mathcal{S}_K^\prime \|_{\varepsilon(\alpha)}^2$ and 
$\| \mathcal{S}_L^\prime \|_{\varepsilon(\alpha)}^2$ are obtained with the same matrices $K$ and $L$, specifically:
\begin{alignat*}{2}
    & \min_K \| \mathcal{S}_K^\prime \|_{\varepsilon(\alpha)}^2 && = \operatorname{trace} (D_1^\top L Q_\alpha L D_1), \\
    & \min_L \| \mathcal{S}_L^\prime \|_{\varepsilon(\alpha)}^2 && = \operatorname{trace} (D_2 K P_\alpha K^\top D_2^\top). 
    \end{alignat*}
The claim of Theorem \ref{thm:KL} follows.
\endproof
\end{prooftheorem}

\bibliographystyle{IEEEtran} 
\bibliography{my}

\begin{thebibliography}{10}
\providecommand{\url}[1]{#1}
\csname url@samestyle\endcsname
\providecommand{\newblock}{\relax}
\providecommand{\bibinfo}[2]{#2}
\providecommand{\BIBentrySTDinterwordspacing}{\spaceskip=0pt\relax}
\providecommand{\BIBentryALTinterwordstretchfactor}{4}
\providecommand{\BIBentryALTinterwordspacing}{\spaceskip=\fontdimen2\font plus
\BIBentryALTinterwordstretchfactor\fontdimen3\font minus
  \fontdimen4\font\relax}
\providecommand{\BIBforeignlanguage}[2]{{%
\expandafter\ifx\csname l@#1\endcsname\relax
\typeout{** WARNING: IEEEtran.bst: No hyphenation pattern has been}%
\typeout{** loaded for the language `#1'. Using the pattern for}%
\typeout{** the default language instead.}%
\else
\language=\csname l@#1\endcsname
\fi
#2}}
\providecommand{\BIBdecl}{\relax}
\BIBdecl

\bibitem{b12-boyd}
S.~{Boyd}, L.~{El Ghaoui}, E.~{Feron}, and V.~{Balakrishnan}, \emph{Linear
  Matrix Inequalities in System and Control Theory}.\hskip 1em plus 0.5em minus
  0.4em\relax SIAM studies in applied mathematics, 1994.

\bibitem{b5-abedor}
J.~Abedor, K.~M. Nagpal, and K.~Poolla, ``A linear matrix inequality approach
  to peak‐to‐peak gain minimization,'' \emph{International Journal of
  Robust and Nonlinear Control}, vol.~6, pp. 899--927, 1996.

\bibitem{b1-nazin}
S.~A. Nazin, B.~T. Polyak, and M.~V. Topunov, ``Rejection of bounded exogenous
  disturbances by the method of invariant ellipsoids,'' \emph{Automation and
  Remote Control}, vol.~68, no.~3, pp. 467--486, 2007.

\bibitem{b13-khlebnikov_filter}
M.~V. Khlebnikov, B.~T. Polyak, and V.~M. Kuntsevich, ``Optimization of linear
  systems subject to bounded exogenous disturbances: The invariant ellipsoid
  technique,'' \emph{Automation and Remote Control}, vol.~72, no.~11, pp.
  2227--2275, 2011.

\bibitem{b2-topunov}
B.~T. Polyak and M.~V. Topunov, ``Suppression of bounded exogenous
  disturbances: Output feedback,'' \emph{Automation and Remote Control},
  vol.~69, no.~5, pp. 801--818, 2008.

\bibitem{e1}
S.~Gonzalez-Garcia, A.~Polyakov, and A.~Poznyak, ``Output linear controller for
  a class of nonlinear systems using the invariant ellipsoid technique,'' in
  \emph{2009 American Control Conference}, 2009, pp. 1160--1165.

\bibitem{x2}
V.~Azhmyakov, A.~Poznyak, and R.~Juárez, ``On the practical stability of
  control processes governed by implicit differential equations: The invariant
  ellipsoid based approach,'' \emph{Journal of the Franklin Institute}, vol.
  350, no.~8, pp. 2229--2243, 2013.

\bibitem{y1}
P.~Ordaz and A.~Poznyak, ``The {F}uruta's pendulum stabilization without the
  use of a mathematical model: Attractive ellipsoid method with
  {KL}-adaptation,'' in \emph{51st IEEE Conference on Decision and Control},
  2012, pp. 7285--7290.

\bibitem{x6}
V.~Azhmyakov, ``On the geometric aspects of the invariant ellipsoid method:
  Application to the robust control design,'' in \emph{50th IEEE Conference on
  Decision and Control}, 2011, pp. 1353--1358.

\bibitem{x7}
V.~Azhmyakov, A.~Poznyak, and O.~Gonzalez, ``On the robust control design for a
  class of nonlinearly affine control systems: The attractive ellipsoid
  approach,'' \emph{Journal of Industrial and Management Optimization}, vol.~9,
  no.~3, pp. 579--593, 2013.

\bibitem{e3}
A.~Polyakov and A.~Poznyak, ``Invariant ellipsoid method for minimization of
  unmatched disturbances effects in sliding mode control,'' \emph{Automatica},
  vol.~47, no.~7, pp. 1450--1454, 2011.

\bibitem{y3}
J.~Davila and A.~Poznyak, ``Dynamic sliding mode control design using
  attracting ellipsoid method,'' \emph{Automatica}, vol.~47, no.~7, pp.
  1467--1472, 2011.

\bibitem{b3-poznyak}
A.~Poznyak, A.~Polyakov, and V.~Azhmyakov, \emph{Attractive Ellipsoids in
  Robust Control}.\hskip 1em plus 0.5em minus 0.4em\relax Springer
  International Publishing, 2014.

\bibitem{y6}
H.~Alazki, E.~Hern{\'a}ndez, J.~M. Ibarra, and A.~Poznyak, ``Attractive
  ellipsoid method controller under noised measurements for {SLAM},''
  \emph{International Journal of Control, Automation and Systems}, vol.~15,
  no.~6, pp. 2764--2775, 2017.

\bibitem{e5}
P.~García and K.~Ampountolas, ``Robust disturbance rejection by the attractive
  ellipsoid method – part {I}: Continuous-time systems,''
  \emph{IFAC-PapersOnLine}, vol.~51, no.~32, pp. 34--39, 2018.

\bibitem{e6}
P.~Garc{í}a and K.~Ampountolas, ``Robust disturbance rejection by the
  attractive ellipsoid method – part {II}: Discrete-time systems,''
  \emph{IFAC-PapersOnLine}, vol.~51, no.~32, pp. 93--98, 2018.

\bibitem{y4}
F.~Oliva-Palomo, A.~Sanchez-Orta, P.~Castillo, and H.~Alazki, ``Nonlinear
  ellipsoid based attitude control for aggressive trajectories in a quadrotor:
  Closed-loop multi-flips implementation,'' \emph{Control Engineering
  Practice}, vol.~77, pp. 150--161, 2018.

\bibitem{y5}
R.~Falcón, H.~Ríos, M.~Mera, and A.~Dzul, ``Attractive ellipsoid-based robust
  control for quadrotor tracking,'' \emph{IEEE Transactions on Industrial
  Electronics}, vol.~67, no.~9, pp. 7851--7860, 2020.

\bibitem{y7}
P.~García and K.~Ampountolas, ``Robust stability of time-varying polytopic
  systems by the attractive ellipsoid method,'' in \emph{American Control
  Conference}, 2019, pp. 1139--1144.

\bibitem{b9-optimalbook}
B.~D. Anderson and J.~B. Moore, \emph{Linear optimal control,
  Prentice-Hall}.\hskip 1em plus 0.5em minus 0.4em\relax Inc., New Jersey,
  1971.

\bibitem{b16-usual_ellipsoids}
G.~Dullerud and F.~Paganini, \emph{A Course in Robust Control Theory: A Convex
  Approach}.\hskip 1em plus 0.5em minus 0.4em\relax Springer New York, 2013.

\bibitem{b17-h2_norm}
R.~Skelton, T.~Iwasaki, and K.~Grigoriadis, \emph{A unified algebraic approach
  to linear control design}.\hskip 1em plus 0.5em minus 0.4em\relax London:
  Taylor \& Francis, 1998.

\bibitem{b6-starcapture}
S.~Venkatesh and M.~Dahleh, ``Does star norm capture $l_1$ norm?'' in
  \emph{Proceedings of 1995 American Control Conference}, vol.~1, 1995, pp.
  944--945.

\bibitem{b7-support1}
A.~Girard and C.~L. Guernic, ``Efficient reachability analysis for linear
  systems using support functions,'' \emph{IFAC Proceedings Volumes}, vol.~41,
  no.~2, pp. 8966--8971, 2008, 17th IFAC World Congress.

\bibitem{b8-support2}
E.~Goncharova and A.~Ovseevich, ``Small-time reachable sets of linear systems
  with integral control constraints: Birth of the shape of a reachable set,''
  \emph{Journal of Optimization Theory and Applications}, vol. 168, 2015.

\bibitem{b10-h2}
J.~Doyle, K.~Glover, P.~Khargonekar, and B.~Francis, ``State-space solutions to
  standard $\mathcal{H}_2$ and $\mathcal{H}_\infty$ control problems,''
  \emph{IEEE Transactions on Automatic Control}, vol.~34, no.~8, pp. 831--847,
  1989.

\bibitem{b11-mathbook}
H.~Trentelman, A.~Stoorvogel, and M.~Hautus, \emph{Control Theory for Linear
  Systems}.\hskip 1em plus 0.5em minus 0.4em\relax Springer London, 2012.

\bibitem{b15-riccati_comparison}
A.~Ran and R.~Vreugdenhil, ``Existence and comparison theorems for algebraic
  {Riccati} equations for continuous- and discrete-time systems,'' \emph{Linear
  Algebra and its Applications}, vol.~99, pp. 63--83, 1988.

\bibitem{b14-riccati_equation}
M.-L. Ni, ``Existence condition on solutions to the algebraic {Riccati}
  equation,'' \emph{Acta Automatica Sinica}, vol.~34, no.~1, pp. 85--87, 2008.

\end{thebibliography}

\vspace{-0.6cm}

\begin{IEEEbiography}[{\includegraphics[width=1in,height=1.25in,clip, keepaspectratio]{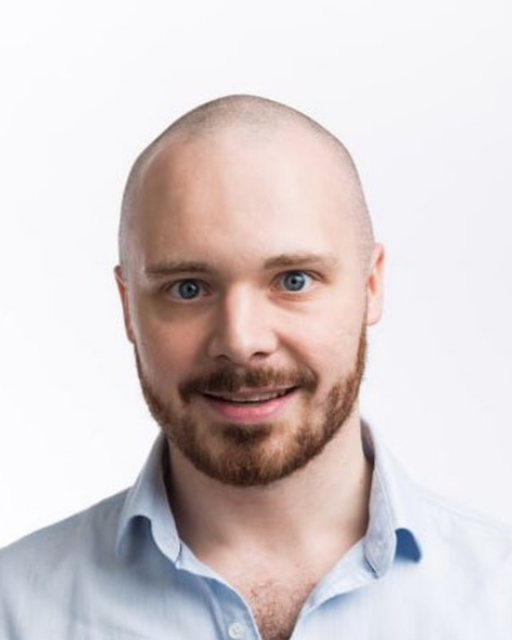}}]{Alexey Peregudin} received the M.Sc. degree in mechanical engineering and the Candidate of Sciences (Ph.D.) degree in control systems engineering from ITMO University, Saint Petersburg, Russia, in 2017 and 2021, respectively.

He is currently a Researcher with the Laboratory of Adaptive and Intelligent Control of Networked and Distributed Systems at the Institute for Problems in Mechanical Engineering, Saint Petersburg, Russia. His research interest includes linear control theory, optimal control, robust control, hybrid systems, finite-time and fixed-time stability, and data-driven control. 

Dr. Peregudin was a recipient of the Best Report Award and the Professor Nesenyuk Prize in 2021 at the Conference on Navigation and Motion Control, and a winner of the St. Petersburg Grant Competition for students, graduate students,
young scientists and candidates of sciences in 2018 and 2021.

\end{IEEEbiography}

\vspace{-0.8cm}

\begin{IEEEbiography}[{\includegraphics[width=1in,height=1.25in,clip, keepaspectratio]{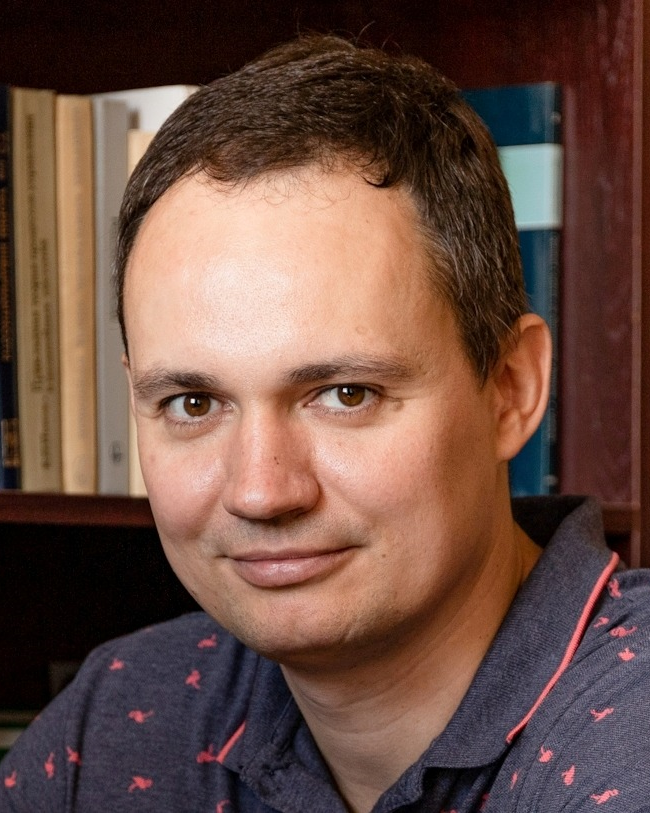}}]{Igor Furtat } 
was born in the former Soviet Union in 1983. He received the Diploma degree in 2005, Candidate of technical sciences (Ph.D.) degree in 2006, Doctor of technical sciences (Habilitation) degree in 2012, Professor in 2018.

He is currently a Head of laboratory in the Institute of problems of mechanical engineering of the Russian academy of sciences. His research interests include nonlinear control, adaptive control, optimal and robust control, systems with time-delays, control of distributed systems, control of dynamical networks, control in the chemical industry, and in power systems.
Since 2015 IEEE member and since 2018 IEEE Senior member.
Since 2017 associate editor of IEEE American Control Conference and IEEE Conference on Decision and Control.

Dr. Furtat was Awarded by the Government of St. Petersburg (Russia) in the field of scientific and pedagogical activity in 2015. Awarded a medal of the Russian Academy of Sciences for development of disturbance compensation theory in 2016.
Best paper award at the IEEE 9th International Congress on Ultra Modern Telecommunications and Control Systems and Workshops, November 6-8, Munich, Germany, 2017.
\end{IEEEbiography}

\end{document}